\documentclass[11pt]{article}

\usepackage{amsmath,amssymb,amsfonts,enumitem,amsthm,accents,pstricks,pst-node,pstricks-add,mathrsfs,xy}
\usepackage{pst-solides3d}

\xyoption{all}

\usepackage{color}
\usepackage{makeidx}
\usepackage{graphicx}
\usepackage{setspace}
\usepackage{extarrows}
\usepackage{hyperref}

\numberwithin{equation}{section}

\usepackage{tocloft}

\usepackage[Symbolsmallscale]{upgreek} 

\newcommand*{\defeq}{\mathrel{\vcenter{\baselineskip0.5ex \lineskiplimit0pt
                     \hbox{\scriptsize.}\hbox{\scriptsize.}}}%
                     =}
                     
\newcommand{\sslash}{\!\mathbin{/\mkern-6mu/\!}}

\newtheorem{theorem}{Theorem}[section]
\newtheorem{lemma}[theorem]{Lemma}
\newtheorem{conjlemma}[theorem]{Conjectural Lemma}
\newtheorem{corollary}[theorem]{Corollary}
\newtheorem{proposition}[theorem]{Proposition}
\newtheorem{conjecture}[theorem]{Conjecture}

\theoremstyle{definition}
\newtheorem{definition}[theorem]{Definition}
\newtheorem{remark}[theorem]{Remark}
\newtheorem{example}[theorem]{Example}

\newtheorem{convention}[theorem]{Convention}

\def\ov#1{\overline{#1}}
\def\tn#1{\textnormal{#1}}
\def\mf#1{\mathfrak{#1}}
\def\wt#1{\widetilde{#1}}

\def\ll{\left\langle}
\def\rr{\right\rangle}
\def\mc{\mathcal}
\def\lra{\longrightarrow}

\newcommand{\abs}[1]{\left\vert #1 \right\vert}
\newcommand{\lrp}[1]{\left( #1 \right)}
\newcommand{\lrc}[1]{\left\{ #1 \right\}}

\def\bEq#1{\begin{equation}\label{#1}}
\def\eEq{\end{equation}}

\def\bsEq{\begin{equation*}}
\def\esEq{\end{equation*}}

\def\bDf#1{\begin{definition}\label{#1}}
\def\eDf{\end{definition}}

\def\bTh#1{\begin{theorem}\label{#1}}
\def\eTh{\end{theorem}}

\def\bCn#1{\begin{conjecture}\label{#1}}
\def\eCn{\end{conjecture}}

\def\bLm#1{\begin{lemma}\label{#1}}
\def\eLm{\end{lemma}}

\def\bCLm#1{\begin{conjlemma}\label{#1}}
\def\eCLm{\end{conjlemma}}

\def\bRm#1{\begin{remark}\label{#1}}
\def\eRm{\end{remark}}

\def\bEx#1{\begin{example}\label{#1}}
\def\eEx{\end{example}}

\def\bPr#1{\begin{proposition}\label{#1}}
\def\ePr{\end{proposition}}

\def\bCr#1{\begin{corollary}\label{#1}}
\def\eCr{\end{corollary}}

\def\bFg#1{\begin{figure}\label{#1}}
\def\eFg{\end{figure}}

\def\bPf{\begin{proof}}
\def\ePf{\end{proof}}

\def\bIt{\begin{itemize}[leftmargin=*]}
\def\eIt{\end{itemize}}

\def\bEn{\begin{enumerate}[label=$(\arabic*)$,leftmargin=*]}
\def\eEn{\end{enumerate}}

\def\bEnalph{\begin{enumerate}[label=$(\alph*)$,leftmargin=*]}
\def\eEnalph{\end{enumerate}}

\def\cC{\mc{C}}
\def\cD{\mc{D}}
\def\cF{\mc{F}}

\def\cM{\mc{M}}

\def\cX{\mc{X}}
\def\cP{\mc{P}}
\def\cQ{\mc{Q}}

\def\cR{\mc{R}}

\def\cS{\mc{S}}

\def\cC{\mc{C}}

\def\R{\mathbb R}
\def\C{\mathbb C}
\def\Z{\mathbb Z}
\def\Q{\mathbb Q}
\def\P{\mathbb P}

\def\N{\mathbb N}

\def\mft{\mf{t}}

\def\la{\lambda}
\def\La{\Lambda}

\def\De{\Delta}
\def\de{\delta}

\def\Si{\Sigma}
\def\al{\alpha}
\def\be{\beta}

\def\G{\ifmmode{\cal G}\else{${\cal G}$}\fi}        
\def\spin{\ifmmode{\mathrm{spin}}\else{spin}\fi}
\def\Spin{\ifmmode{\mathrm{Spin}}\else{Spin}\fi}
\def\spinc{\ifmmode{\mathrm{spin}^c}\else{spin$^c$ }\fi}
\def\Spinc{\ifmmode{\mathrm{Spin}^c}\else{Spin$^c$ }\fi}

\def\sfrak{\ifmmode{\mathfrak s}\else{${\mathfrak s }$}\fi}  
\def\tfrak{\ifmmode{\mathfrak t}\else{${\mathfrak t }$}\fi} 
\def\sz{\ifmmode{{\tn\ss}}\else{\ss}\fi} 
\def\sbold{\ifmmode{\mbox{\boldmath$s$\unboldmath}}\else{\boldmath$s$\unboldmath}\fi}
\def\tbold{\ifmmode{\mbox{\boldmath$t$\unboldmath}}\else{\boldmath$t$\unboldmath}\fi} 

\usepackage{pict2e}

\makeatletter

\newcommand*\@KP@Large@frame[2]{%
    \setlength\unitlength{\fontdimen 22 #1\tw@}%
    \vrule \@width\z@ \@height 4\unitlength \@depth\tw@\unitlength
    \begin{picture}(6,2)(-3,-1)%
        \def\@KP@Radius     {3}%
        \def\@KP@Hole@radius{.5}
        \def\@KP@Diameter   {6}%
        #2%
    \end{picture}%
}
\newcommand*\@KP@Small@frame[2]{%
    \setlength\unitlength{\fontdimen 22 #1\tw@}%
    \vrule \@width\z@ \@height \thr@@\unitlength \@depth\@ne\unitlength
    \begin{picture}(4,2)(-2,-1)%
        \def\@KP@Radius     {2}%
        \def\@KP@Hole@radius{.5}
        \def\@KP@Diameter   {4}%
        #2%
    \end{picture}%
}

\newcommand*\@KP@Radius     {}
\newcommand*\@KP@Hole@radius{}
\newcommand*\@KP@Diameter   {}
\newcommand*\@KP@Shape@A{%
    \put(0,0){\circle{\@KP@Diameter}}%
}
\newcommand*\@KP@Shape@B{%
    \Line(-\@KP@Radius,\@KP@Radius )(\@KP@Radius,-\@KP@Radius)%
    \Line(-\@KP@Radius,-\@KP@Radius)(-\@KP@Hole@radius,-\@KP@Hole@radius)%
    \Line(\@KP@Radius ,\@KP@Radius )(\@KP@Hole@radius ,\@KP@Hole@radius )%
}
\newcommand*\@KP@Shape@E{%
    \Line(-\@KP@Radius,\@KP@Radius )(\@KP@Radius,-\@KP@Radius)%
     \Line(\@KP@Radius,\@KP@Radius )(-\@KP@Radius,-\@KP@Radius)%

}

\newcommand*\@KP@Shape@C{%
    \cbezier(-\@KP@Radius,\@KP@Radius )(0,0)(0,0)(\@KP@Radius,\@KP@Radius )%
    \cbezier(-\@KP@Radius,-\@KP@Radius)(0,0)(0,0)(\@KP@Radius,-\@KP@Radius)%
}
\newcommand*\@KP@Shape@D{%
    \cbezier(-\@KP@Radius,-\@KP@Radius)(0,0)(0,0)(-\@KP@Radius,\@KP@Radius)%
    \cbezier(\@KP@Radius ,-\@KP@Radius)(0,0)(0,0)(\@KP@Radius ,\@KP@Radius)%
}

\newcommand*\@KP@Atomic@mathpalette[1]{%
    \mathinner{
        \mathchoice{%
            \linethickness{.6\p@}
            \@KP@Large@frame \textfont {#1}%
        }{%
            \linethickness{.4\p@}
            \@KP@Small@frame \textfont {#1}%
        }{%
            \linethickness{.3\p@}
            \@KP@Small@frame \scriptfont {#1}%
        }{%
            \linethickness{.2\p@}
            \@KP@Small@frame \scriptscriptfont {#1}%
        }%
    }%
}

\newcommand*\KPA{\@KP@Atomic@mathpalette \@KP@Shape@A}
\newcommand*\KPB{\@KP@Atomic@mathpalette \@KP@Shape@B}
\newcommand*\KPC{\@KP@Atomic@mathpalette \@KP@Shape@C}
\newcommand*\KPD{\@KP@Atomic@mathpalette \@KP@Shape@D}
\newcommand*\KPE{\@KP@Atomic@mathpalette \@KP@Shape@E}

\makeatother

\makeindex
\begin{document}
\title{The Geometric P=W conjecture \\
and Thurston's compactification}
\author{Ashwin Ayilliath-Kutteri, Mohammad Farajzadeh-Tehrani, Charles Frohman}
\date{\today}
\maketitle

\begin{abstract}
In this paper, we use new results together with established facts about Thurston's compactification of Teichm\"uller space to address the geometric P=W conjecture for $\mathrm{SL}(2,\mathbb{C})$, which concerns projective compactifications of character varieties of closed surfaces. In particular, we construct a projective compactification of the $\mathrm{SL}(2,\mathbb{C})$-character variety of any closed surface of genus $g>1$, in which the boundary divisors are toric varieties and the dual intersection complex is a sphere. A main technical step, of independent interest, is the derivation of an explicit formula for a well-known embedding of the set of isotopy classes of multicurves on a closed surface of genus $g$ into $\mathbb{N}^{9g-9}$.
\end{abstract}

\tableofcontents

\section{Introduction and the main results}
For any reductive Lie group $G$ and finitely generated group $\Gamma$, the $G$-\textbf{character variety}
$$
\cX(\Gamma,G)=\tn{Hom}(\Gamma, G)\sslash G
$$
is the moduli space of $G$-representations of $\Gamma$ up to conjugation by elements of $G$. The special case $\cX_{g}(G)$, where $\Gamma=\pi_1(\Si_{g})$ is the fundamental group of a genus $g$ surface, is well studied and connects to several areas of mathematics. For instance, let $\cM_{g}(G)$ denote the so-called \textbf{Dolbeault moduli space}, whose points parametrize flat principal $G$-Higgs bundles on a closed Riemann surface of genus~$g$. The Dolbeault moduli space is hyperK\"ahler and admits a canonical holomorphic (projective) Lagrangian fibration onto an affine space of half the dimension. There exists a canonical real-analytic identification
$$
\Psi\colon \cM_{g}(G)\lra \cX_{g}(G),
$$
called the \emph{non-abelian Hodge correspondence}.\medskip

The geometric P=W conjecture, proposed by Katzarkov et al.~\cite[Conj.~1.1]{KNPS} and further refined in~\cite{MMS}, aims to provide a deeper understanding of the asymptotic behavior of the transcendental map $\Psi$ with respect to suitable compactifications on both sides.

\begin{conjecture}[\cite{MMS}]\label{PWGeometric}
$\cX_{g}(G)$ admits a dlt log Calabi--Yau compactification, and the dual intersection complex of the boundary divisor is a polyhedral complex homeomorphic to a sphere.
\end{conjecture}

\begin{remark}
The more sophisticated homological P=W conjecture~\cite{DHM} asserts that the non-abelian Hodge correspondence exchanges the perverse and weight filtrations. Under certain technical assumptions, the geometric P=W conjecture implies the cohomological P=W conjecture at the highest weight level. It is also worth noting that Katzarkov et al.'s original conjecture asserts that the dual intersection complex of the boundary divisor is homotopy equivalent to a sphere. This weaker version has recently been proved in the punctured case for all $G$; cf.~\cite{Su}.
\end{remark}

Conjecture~\ref{PWGeometric} extends to punctured Riemann surfaces and involves the interplay between ``relative character varieties" and parabolic Higgs bundles. In recent work~\cite{FF23}, Frohman and Tehrani used
\begin{itemize}
\item the canonical isomorphism between the ring of regular functions on the $\tn{SL}(2,\C)$ character variety of a punctured surface and its skein algebra, and
\item a filtration on the skein algebra arising from intersection numbers with the arcs of an ideal triangulation
\end{itemize}
to construct compactifications of character varieties of punctured surfaces with interesting properties.

\medskip

Let $\Si_{g,n}$ denote an $n$-punctured genus $g$ surface. For $G=\tn{SL}(2,\C)$, the negative of the trace function along curves around the punctures defines a fibration
\begin{equation}\label{Rel-family}
\pi\colon  \cX_{g,n}:= \cX\big(\pi_{1}(\Si_{g,n}),\tn{SL}_2(\C)\big) \lra \C^n, 
\end{equation}
whose fibers $\cX_{g,n,\mf{t}}=\pi^{-1}(\mf{t})$ are called the \textbf{relative $\tn{SL}(2,\C)$-character varieties}. It is well known that each $\cX_{g,n,\mf{t}}$ is quasi-projective, hyperK\"ahler, and log Calabi--Yau (though, sometimes singular). For other choices of $G$, the relative character varieties are often defined by fixing the conjugacy classes of the representations at the punctures. 

\begin{theorem}[\cite{FF23}]\label{FF23theorem}
For $g\geq 0$ and $n>0$ satisfying $2g+n\geq 3$, associated to any ``ideal triangulation" of $\Si_{g,n}$, there is a normal compactification $\ov{\cX}_{g,n}$ of $\cX_{g,n}$ such that the boundary divisor $D_{g,n}=\partial \ov{\cX}_{g,n}$ is an irreducible ample toric variety. Additionally, there is a normal relative (i.e.\ fiberwise) compactification $\ov{\cX}^{\tn{rel}}_{g,n}$ of the family~\eqref{Rel-family} such that, for each $\mft\in \C^n$, the boundary divisor $D_{g,n,\mft}$ is a toric and $\mf{t}$-independent subvariety of $D_{g,n}$. Furthermore, the moment polytope complex of $D_{g,n,\mft}$ is a sphere.
\end{theorem}

Here, we give a quick tour of the proof of Theorem~\ref{FF23theorem} which combines several major results from quantum topology. 

\medskip

Let $\mathcal{S}_{g,n}$ denote the set of isotopy classes of nontrivial\footnote{i.e., curves that do not bound a disk.} simple closed curves (embedded circles) on $\Sigma_{g,n}$. For $c_1, c_2 \in \mathcal{S}_{g,n}$, let $i(c_1, c_2)$ denote their geometric intersection number, defined as the minimum number of intersection points between a simple curve $c'_1$ isotopic to $c_1$ and a simple curve $c'_2$ isotopic to $c_2$, where $c'_1$ is transverse to $c'_2$. More generally, let $\mathcal{S}'_{g,n}$ denote the set of isotopy classes of simple, closed, but not necessarily connected curves on $\Sigma_{g,n}$, called \textbf{multicurves}\footnote{also called ``simple diagrams''.}, whose components all belong to $\mathcal{S}_{g,n}$. The notion of geometric intersection number naturally extends to a symmetric pairing on $\mathcal{S}'_{g,n}$.

\medskip

As we explain in more detail in Section~\ref{sec:SK-DT}, by a theorem of Bullock--Przytycki--Sikora~\cite{Bu,PS} and Charles--March\'{e}~\cite{ChMa}, via the (negative of the) trace function along loops, the ring of regular functions $\C[\cX_{g,n}]$ on $\cX_{g,n}$ is canonically isomorphic to the classical skein algebra $\tn{Sk}_{g,n}$ of the surface. Furthermore, by~\cite[Thm.~IX.7.1]{P1}, $\mathcal{S}'_{g,n}$ forms a vector space basis of  $\tn{Sk}_{g,n}$, and (if $n>0$) taking intersection numbers with the arcs of an ideal triangulation defines an embedding
$$
\iota \colon \mathcal{S}'_{g,n}\to \mathbb{N}^{6g-6+3n}
$$
whose image lies in an explicitly characterized rational polyhedral cone $\Lambda$. This embedding defines a projective filtration on $\C[\cX_{g,n}]$ and induces an isomorphism between the graded coordinate ring of $\cX_{g,n}$ and the algebra $\C[\iota(\mathcal{S}'_{g,n})]$ associated to the (saturated) monoid $\iota(\mathcal{S}'_{g,n}) \subset \Lambda$. It follows that the projective compactification obtained from this filtration adds an irreducible divisor $D_{g,n}$ that is a projective toric variety with moment polytope $P$, such that $\La$ is a cone over $P$.

\medskip

\emph{In this paper, using intersection numbers with closed curves, we prove an analogous result for closed Riemann surfaces. The proof, however, is significantly more involved and uses a stronger version of the famous $(9g-9)$-Theorem stated below that yields an embedding map
$$
\iota \colon \mathcal{S}'_g \to \mathbb{N}^{9g-9}
$$ }
with similar properties.

\medskip

Closely related to set of multicurves $\cS'_g$ on $\Si_g$ is the space\footnote{Denoted by $\cM\cF$ in~\cite{FLP}.} $\mathcal{F}_g$ of (Whitehead equivalence classes of) \textbf{measured foliations} on $\Sigma_g$. This space admits a natural inclusion $\mathcal{S}'_g \subset \mathcal{F}_g$ such that taking intersection numbers with curves in $\mathcal{S}_g$ extends continuously to $\mathcal{F}_g$. In other words, elements of $\mathcal{F}_g$ may be viewed as multicurves endowed with positive real weights.

\medskip

A main technical step of the independent interest to topologists is an explicit version of the well-known $(9g-9)$-Theorem; see~\cite[Theorem~4.10]{FLP} or~\cite[p.~438]{FM}, which can be briefly stated as follows.

\begin{theorem}\label{ConeTheorem}
For $g>1$, there exists a collection of $9g-9$ curves $\mc{C}$ on $\Si_g$ such that taking intersection numbers with the curves in $\mc{C}$ defines an embedding $\iota_\cC \colon \mathcal{S}'_g \to \mathbb{N}^{9g-9}$ with the following properties.
\begin{itemize}
\item The embedding above extends continuously to an embedding
$$
\iota_\cC \colon \mathcal{F}_g\lra \R_{\geq 0}^{9g-9}
$$
whose image is a union $\Lambda=\bigcup_{\ell}\La_\ell$ of $(6g-6)$-dimensional convex rational polyhedral cones over a polyhedral complex
$$
P=\bigcup_{\ell} P_\ell
$$
homeomorphic to $S^{6g-7}$. 
\item The image of $\cS'_g$ in each $\Lambda_\ell$ is a finite-index (saturated) submonoid $Q_\ell$ of $\Lambda_\ell\cap \N^{9g-9}$.
\item Restricted to the preimage of each $\Lambda_\ell$, $\iota_\cC$ is linear in Dehn--Thurston (as well as ``modified" Dehn--Thurston) coordinates.
\end{itemize}
\end{theorem}

Then, by comparing the multiplication of multicurves in the skein algebra $\tn{Sk}_g$ with the additive monoid structure of each $Q_\ell$, we identify the graded algebra of $\C[\cX_{g}]$ with $\C[\iota_\cC(\cS'_g)]$ and prove the following result.

\begin{theorem}\label{MainTheorem}
Let $\mc{C}$ be the collection of $9g-9$ curves in the (proof) of Theorem~\ref{ConeTheorem}.
For every multicurve $\mf{m}$ on $\Si_g$, define its degree by
\begin{equation}\label{grading0}
|\mf{m}|\defeq \sum_{c\in \cC}  i(\mf{m},c).
\end{equation}
Under the canonical isomorphism between the ring of regular functions $\C[\cX_{g}]$ and the classical skein algebra of the surface $\tn{Sk}_g$, the degree function~\eqref{grading0} defines a projective filtration on $\C[\cX_{g}]$ and thus a normal projective compactification $\ov{\cX}_{g}$ of $\cX_{g}$ such that the boundary divisor $D_g=\partial \ov{\cX}_{g}$ is a union of toric varieties intersecting along toric strata according to the polytope complex $P$  in Theorem~\ref{ConeTheorem}. Consequently, the intersection complex of $D_{g}$ is dual to the moment polytope complex $P$ and  is also homeomorphic to $S^{6g-7}$.
\end{theorem}

\begin{remark}
In general, the resulting compactification $\ov\cX_g$ need not be dlt, but it can be made so by performing (toric) blowups that preserve the toric structure of the boundary divisor and maintain the polytope complex as a sphere. The log Calabi--Yau (log CY) property, however, is only known for punctured surfaces (see~\cite{Wang}) and remains an open question for closed surfaces.
\end{remark}

\textbf{Outline of the paper.} In Section~\ref{sec:SK-DT}, we give a brief introduction to the definition of the skein algebra and Dehn--Thurston coordinates. Readers already familiar with these notions may skip ahead to Section~\ref{sec:TC}. In Section~\ref{sec:TC}, we review the construction of Thurston's compactification of Teichm\"uller space and   state a theorem that provides an explicit formula for intersection numbers with the curves in a distinguished collection of $9g-9$ curves in terms of Dehn--Thurston coordinates. We then use this formula to prove Theorem~\ref{ConeTheorem}. We also explain the connection of Conjecture~\ref{PWGeometric} with Thurston's compactification and provide remarks on the significance of our formulas. In Section~\ref{sec:conenction}, we describe the filtration arising from intersection numbers, briefly review the connection between monoids and toric varieties, and prove the main result--Theorem~\ref{MainTheorem}. Sections~\ref{sec:DTandf-functions} and~\ref{Sec:proof-of-MainTheorem} contain the proofs of the technical statements concerning intersection coordinates and their behavior under the product map in the skein algebra. Finally, in Section~\ref{sec:g2}, we discuss the example of genus two and give explicit description of toric divisors.

\medskip 

{\bf Acknowledgement}.We are grateful to Wade Bloomquist for informing us about the results of~\cite{ABCGLW}, which were inspiring. We also thank Robert Lipshitz for sharing some anonymous comments with us that were helpful in revising the manuscript.

\section{Preliminaries in skein algebra and DT coordinates}\label{sec:SK-DT}

In this section, we give a brief introduction to the definition of the skein algebra and Dehn--Thurston coordinates, focusing on closed surfaces for simplicity. \\

\textbf{Skein algebra.} For $g\geq 0$, the (quantum) \textbf{skein algebra} $\tn{Sk}_{q}(\Si_{g})$ of $\Si_{g}$ is the $\Z[q^{\pm}]$-algebra generated by isotopy classes of framed links in $\Si_{g}\times \R$, modulo the relations
\[
\aligned \left<L\cup\KPA\right> &= (-q^{2}-q^{-2})\langle L\rangle &&\text{disjoint union with trivial loop} \\ \left<\KPB\right> &= q\left<\KPC\right> + q^{-1}\left<\KPD\right> &&\text{skein relation,} \endaligned
\]
with the product structure defined by stacking two links on top of each other. By~\cite[Thm.~IX.7.1]{P1}, isotopy classes of multicurves $\cS'_g$ (together with the empty multicurve, identified with $1$) form a basis of $\tn{Sk}_{q}(\Si_{g})$ as a $\Z[q^{\pm}]$-module.

On the other hand, by a theorem of Bullock--Przytycki--Sikora~\cite{Bu,PS} and Charles--March\'{e}~\cite{ChMa}, the specialization
\[
\tn{Sk}_{g}\defeq \tn{Sk}_{q=-1}(\Si_{g})
\]
of the skein algebra at $q=-1$ is canonically isomorphic, via the (negative of the) trace function along loops, to the ring of regular functions $\C[\cX_{g}]$ on $\cX_{g}$. The classical limit $\tn{Sk}_{q=1}(\Si_{g})$ is isomorphic to $\tn{Sk}_{q=-1}(\Si_{g})$, though not canonically: the isomorphism depends on a choice of spin structure; see~\cite{Barrett}. For simplicity, we work with the $q=1$ classical limit.

In this classical limit, the commutative product structure, denoted by $\ast$, is given by taking the union of two multicurves (in transverse position) and using the skein relations as reduction rules
\begin{equation}\label{skein-relation}
   \left<L\cup\KPA\right> =-2\langle L\rangle,\qquad \left<\KPE\right> = \left<\KPC\right> + \left<\KPD\right>
\end{equation}
to expand the product as a linear combination of multicurves. Therefore, for any two multicurves $\mf{m}$ and $\mf{m}'$, we have
\[
\mf{m}\ast\mf{m}'=\sum_{s} \la_s \mf{m}_s,
\]
where $s$ runs over all possible smoothings of $\mf{m}\cup \mf{m}'$ (in general position), and the coefficient $\la_s$ records the contributions of trivial components (with a trivial curve contributing $-2$) appearing in a smoothing. In later sections, we introduce a filtration on $\tn{Sk}_g$ and give an explicit description of its graded algebra.\\

\textbf{Dehn-Thurston coordinates.} 
Dehn-Thurston (DT, for short) coordinates provide a way to describe elements of $\cS'_g$ using a collection of $6g-6$ integers. Moreover, they allow one to draw standard representatives of multicurves and to use these representatives to study their product in the skein algebra. Here, we review the precise definition of DT coordinates, which uses a pants decomposition $\cP$ of $\Si_g$ together with an embedded dual graph~$\Gamma$. We will need these details in the proofs of our results.
\\

A {\bf pants decomposition} of a closed surface $\Sigma_g$ of genus $g$ is a collection 
$$
\cP =\{a_1, \ldots, a_{3g-3}\}
$$ 
of disjoint, non-trivial curves such that no two are isotopic. This collection cuts $\Sigma_g$ into $2g-2$ pair of pants (i.e., $3$-holed spheres). We allow a pants decomposition to include {\bf folded pants}, where two of the boundary curves have the same image in $\Sigma_g$; see Figure~\ref{folded-pants_fig}. We denote a pants bounded by $a_i,a_j,a_k$ as $F(a_i,a_j,a_k)$. For a folded pants where $a_j=a_k$, we may write $F(a_i,2a_j)$ for clarity.\\

\begin{figure}
\begin{pspicture}(-15,-.8)(0,.8)
\psset{unit=1cm}

\psarc(-7,0){1}{90}{270}\psline(-7,1)(-6,1)\psline(-7,-1)(-6,-1)
\psarc(-7.5,0){.5}{-45}{45}
\psarc(-6.7,0){.5}{145}{215}
\psellipticarc[linecolor=red, linestyle=dashed, dash=1pt](-7.6,0)(0.4,0.1){0}{180} \rput(-7.6,-.3){\small $a$}
\psellipticarc[linecolor=red](-7.6,0)(0.4,0.1){-180}{0}
\psellipse[linecolor=red](-6,0)(0.2,1)

\end{pspicture}
\caption{A folded pair of pants that includes two copies of the curve $a$ on its boundary.}
\label{folded-pants_fig}
\end{figure}

Given a pants decomposition $\cP=\{a_1,\ldots,a_{3g-3}\}$, a {\bf dual graph} $\Gamma$ to $\cP$ is a connected, trivalent graph embedded\footnote{Abstractly, the dual graph is unique. There are, however, different ways for embedding that into $\Sigma_g$.} into $\Sigma_g$. It has exactly $2g-2$ vertices, one in the interior of each pair of pants, and $3g-3$ edges $\{e_1,\ldots,e_{3g-3}\}$ such that $e_i$ intersects $a_i$ transversally at a single point and is disjoint from $a_j$ for $j\neq i$. The pair $(\cP,\Gamma)$ is called a {\bf coordinate datum} for $\Sigma_g$; c.f. \cite[Section~3.6]{FKL}. Note that an edge $e_j$ dual to a curve $a_j$ that is part of a folded pair of pants $F(a_i,2a_j)$ will be a loop in $\Gamma$.\\

Let $N(\Gamma)$ be a regular closed neighborhood (i.e., a thickening) of the dual graph $\Gamma$. This means $N(\Gamma)$ is a subsurface of $\Sigma_g$ with a boundary, containing $\Gamma$ in its interior, such that $\Gamma$ is a strong deformation retract of $N(\Gamma)$. The collection of boundary curves
$$
\partial N(\Gamma)=\{b_1,\ldots,b_\ell\}
$$
defines a decomposition of each pair of pants $F$ into two hexagons. We define the {\bf $\Gamma$-hexagon} as the one belonging to $N(\Gamma)$. Figure~\ref{coordinate-datum_fig} illustrates two pants decompositions of a genus two surface and the curves $b_i$. We always draw the pictures such that the dual graph and thus the $\Gamma$-hexagons are in front. Also, sometimes, we only draw the chosen half.\\

\begin{figure}
\begin{pspicture}(-10.5,-1)(0,1.5)
\psset{unit=1cm}

\psarc[linecolor=blue](-7,0){1}{90}{270}\psline[linecolor=blue](-7,1)(-5,1)\psline[linecolor=blue](-7,-1)(-5,-1)\psarc[linecolor=blue](-5,0){1}{-90}{90}\rput(-6.2,1.2){\small $b_2$}

\psellipse[linecolor=blue](-7.1,0)(.25,.5)\rput(-7.1,.7){\small $b_1$}
\psellipse[linecolor=blue](-4.9,0)(.25,.5)\rput(-4.9,.7){\small $b_3$}

\psarc(-7.5,0){.5}{-45}{45}   \psarc(-4.5,0){.5}{135}{225}
\psarc(-6.7,0){.5}{145}{215} \psarc(-5.3,0){.5}{-35}{35}
\psellipticarc[linecolor=red, linestyle=dashed, dash=1pt](-7.6,0)(0.4,0.1){0}{180} \rput(-7.6,.2){\small $a_1$}
\psellipticarc[linecolor=red](-7.6,0)(0.4,0.1){-180}{0}

\psellipticarc[linecolor=red, linestyle=dashed, dash=1pt](-4.4,0)(0.4,0.1){0}{180} \rput(-4.4,.2){\small $a_3$}
\psellipticarc[linecolor=red](-4.4,0)(0.4,0.1){-180}{0}

\psellipticarc[linecolor=red](-6,0)(0.2,1){90}{270} \rput(-6.4,0){\small $a_2$}
\psellipticarc[linecolor=red, linestyle=dashed, dash=1pt](-6,0)(0.2,1){-90}{90}

\psarc[linecolor=blue](0,0){1}{90}{270}\psline[linecolor=blue](0,1)(2,1)\psline[linecolor=blue](0,-1)(2,-1)\psarc[linecolor=blue](2,0){1}{-90}{90}

\rput(.8,1.2){\small $b_2$}

\psellipse[linecolor=blue](-0.1,0)(.25,.5)\rput(-0.1,.7){\small $b_1$}
\psellipse[linecolor=blue](2.1,0)(.25,.5)\rput(2.1,.7){\small $b_3$}

\psarc(-0.5,0){.5}{-45}{45}   \psarc(2.5,0){.5}{135}{225}
\psarc(.3,0){.5}{145}{215} \psarc(1.7,0){.5}{-35}{35}

\psellipticarc[linecolor=red, linestyle=dashed, dash=1pt](-0.6,0)(0.4,0.1){0}{180} \rput(-0.6,.2){\small $a_1$}
\psellipticarc[linecolor=red](-0.6,0)(0.4,0.1){-180}{0}

\psellipticarc[linecolor=red, linestyle=dashed, dash=1pt](2.6,0)(0.4,0.1){0}{180} \rput(2.6,.2){\small $a_3$}
\psellipticarc[linecolor=red](2.6,0)(0.4,0.1){-180}{0}

\psellipticarc[linecolor=red, linestyle=dashed, dash=1pt](1,0)(1,0.1){0}{180} \rput(1,0.2){\small $a_2$}
\psellipticarc[linecolor=red](1,0)(1,0.1){180}{360}

\end{pspicture}
\caption{Two pants decompositions of a genus two surface.}
\label{coordinate-datum_fig}
\end{figure}

If $F=F(a_1,a_2,a_3)$ is a pants with (ordered) boundary curves $a_1,a_2,a_3$, then
\begin{enumerate}
    \item every simple closed curve in $F$ is homotopic to one of its boundary curves,
    \item and, any mutually disjoint collection $\gamma$ of arcs (a {\bf multiarc}) on $F$ -- up to isotopy -- is classified by its intersection numbers with $a_1,a_2,a_3$:
    \begin{equation}\label{Pants-coordinates_e}
    (i(\gamma,a_1),i(\gamma,a_2),i(\gamma,a_3))\in \lrc{(n_1,n_2,n_3)\in \mathbb{N}^3\colon 2\!\mid\! n_1+n_2+n_3}.
    \end{equation}
\end{enumerate}
Here, we are considering arcs that start and end on a boundary component.
Given a decomposition of $F(a_1,a_2,a_3)$ into two hexagons, with $H$ denoting the selected $\Gamma$-hexagon, for every $(n_1,n_2,n_3)$ such that $2\!\mid\! n_1+n_2+n_3$, the standard configuration of $\gamma$ in $F$ is drawn as follows: 

Assume $a_1,a_2,a_3$ are numbered counter-clockwise with respect to the orientation of $H$. We denote the half of each $a_i$ that belongs to $H$ by $\alpha_i$. We also label the other three edges of $H$ by $\beta_1$, $\beta_2$, and $\beta_3$ such that $\beta_i$ is the opposite edge of $\alpha_i$; see Figure~\ref{Triangular_fig}. Each $\beta_i$ is a segment of some $b_j$. The arcs $\beta_1$, $\beta_2$, and $\beta_3$ are also called {\bf seams} of $F$.\\

Place $n_i$ points on the selected half $\alpha_i$ of $a_i$.
\begin{enumerate}
    \item If $(n_1,n_2,n_3)$ satisfies the triangle inequalities
    \begin{equation}\label{T-ineq_e}
    n_i \leq n_j+n_k \qquad \forall~i=1,2,3, \{j,k\}=\{1,2,3\}-\{i\}
    \end{equation}
    then there exists a unique way to connect those points inside $H$ using parallel copies of seam arcs $\beta_i$ to form a multiarc $\gamma$ whose geometric intersection numbers with the boundary components of $F$ are $(n_1,n_2,n_3)$; see Figure~\ref{Triangular_fig}. We call such arcs the {\bf $\beta$-arcs}, with $\be_i$-arcs being specifically those that are parallel to $\beta_i$.

    \begin{figure}
    \begin{pspicture}(-14,-1.5)(0,1.5)
    \psset{unit=1cm}

    \psline(-7.4,1)(-6.6,1) \rput(-7,1.3){\small $\al_3$}
    \psline(-5.4,1)(-4.6,1) \rput(-5,1.3){\small $\al_2$}
    \psline(-6.4,-1)(-5.6,-1)\rput(-6,-1.3){\small $\al_1$}

    \psarc(-6,1){.6}{180}{360}\rput(-6,.7){\small $\beta_1$}
    \psbezier(-4.6,1)(-4.6,-.5)(-5.6,.5)(-5.6,-1) \rput(-5,-.2){\small $\beta_3$}
    \psbezier(-7.4,1)(-7.4,-.5)(-6.4,.5)(-6.4,-1) \rput(-7,-.2){\small $\beta_2$}

    \psbezier[linecolor=red](-4.8,1)(-4.8,-.5)(-5.8,.5)(-5.8,-1)
    \psbezier[linecolor=red](-4.9,1)(-4.9,-.5)(-5.9,.5)(-5.9,-1)
    \psbezier[linecolor=red](-7.2,1)(-7.2,-.5)(-6.2,.5)(-6.2,-1)
    \psbezier[linecolor=red](-7.1,1)(-7.1,-.5)(-6.1,.5)(-6.1,-1)
    \psarc[linecolor=red](-6,1){.8}{180}{360}

    \end{pspicture}
    \caption{A multi-arc with $(n_1,n_2,n_3)=(4,3,3)$ in the red hexagon.}
    \label{Triangular_fig}
    \end{figure}

    \item If one of the triangle inequalities is not satisfied, say $n_1>n_2+n_3$, we still draw all arcs starting and ending from the selected points on $\alpha_i$. However, $(n_1-(n_2+n_3))/2$ of the arcs in $\gamma$ must leave and return to $H$ at two points on the seam curves. These are called {\bf $U$-arcs}, with $U_i$-arcs being specifically those that start and end on $\alpha_i$. We make the choice so that a $U_i$-arc would go around $a_{i+1}$ (with indices being mod $3$). For instance, a $U_1$-arc starting and ending on $\alpha_1$ will intersect $\beta_1$ and then $\beta_3$; see Figure~\ref{Dominant-xi_fig}. 

    \begin{figure}
    \begin{pspicture}(-14,-1.5)(0,1.5)
    \psset{unit=1cm}

    \psellipse(-7,1)(0.4,0.1) \rput(-7,1.3){\small $a_3$}
    \psellipse(-5,1)(0.4,0.1) \rput(-5,1.3){\small $a_2$}
    \psellipticarc[linestyle=dashed, dash=1pt](-6,-1)(.4,.1){0}{180}\rput(-6,-1.3){\small $a_1$}
    \psellipticarc(-6,-1)(.4,.1){180}{360}\rput(-6,-1.3)

    \psarc(-6,1){.6}{180}{360}\rput(-6,.7){\small $\beta_1$}
    \psbezier(-4.6,1)(-4.6,-.5)(-5.6,.5)(-5.6,-1) \rput(-5,-.2){\small $\beta_3$}
    \psbezier(-7.4,1)(-7.4,-.5)(-6.4,.5)(-6.4,-1) \rput(-7,-.2){\small $\beta_2$}

    \psbezier[linecolor=red](-4.9,.9)(-4.9,-.5)(-5.9,.5)(-5.9,-1.1)
    \psbezier[linecolor=red](-7.1,.9)(-7.1,-.5)(-6.1,.5)(-6.1,-1.1)
    \pscurve[linecolor=red](-6,-1.1)(-6.1,0)(-6,.4)
    \pscurve[linecolor=red,linestyle=dashed, dash=1pt](-6,.4)(-5.6,0.1)(-5.5,-.4)
    \pscurve[linecolor=red](-5.5,-.4)(-5.7,-.7)(-5.8,-1.05)

    \end{pspicture}
    \caption{The standard presentation of a multi-arc with $(n_1,n_2,n_3)=(4,1,1)$.}
    \label{Dominant-xi_fig}
    \end{figure}
\end{enumerate}

With the local pictures as described above in every pants, given a coordinate datum $(\cP,\Gamma)$, the {\bf standard} (up to isotopy) representative of any multicurve $\mathfrak{m}$  is determined by the intersection numbers
$$
n=(n_1,\ldots,n_{3g-3})=\big(i(\mathfrak{m},a_1),\ldots,i(\mathfrak{m},a_{3g-3})\big)\in \mathbb{N}^{3g-3}
$$
and some {\bf twist coordinates}
$$
t=(t_1,\ldots,t_{3g-3})\in \mathbb{Z}^{3g-3}~~\tn{s.t.}~~t_i\geq 0~\tn{whenever}~n_i=0.
$$
The tuple $(n,t)$ is called the {\bf Dehn-Thurston coordinates} of $\mathfrak{m}$, and $\mathfrak{m}(n,t)$ is drawn as follows:

Draw two parallel copies of each $a_i$ to create an annulus thickening $A_i$ of $a_i$. The complement
\bEq{shrunken_eq}
\Sigma'_g=\Sigma_g\setminus \bigcup_{i=1}^{3g-3} A_i
\eEq
is a union of shrunken pairs of pants, each divided into two hexagons by the curves in $\partial N(\Gamma)$; see Figure~\ref{shrunkenPP_fig} that corresponds to Figure~\ref{coordinate-datum_fig}.Right.\\
\begin{figure}
\begin{pspicture}(-7,-1)(0,1.5)
\psset{unit=1cm}

\psarc[linecolor=blue](0,0){1}{90}{270}\psline[linecolor=blue](0,1)(2,1)\psline[linecolor=blue](0,-1)(2,-1)\psarc[linecolor=blue](2,0){1}{-90}{90}

\psellipse[linecolor=blue](-0.1,0)(.15,.5)
\psellipse[linecolor=blue](2.1,0)(.15,.5)

\psellipticarc[linecolor=red, linestyle=dashed, dash=1pt](-0.6,0.2)(0.4,0.1){0}{180} \psellipticarc[linecolor=red](-0.6,0.2)(0.4,0.1){-180}{0}

\rput(-.6,0){\small $A_1$}

\psellipticarc[linecolor=red, linestyle=dashed, dash=1pt](-0.6,-0.2)(0.4,0.1){0}{180} \psellipticarc[linecolor=red](-0.6,-0.2)(0.4,0.1){-180}{0}

\psellipticarc[linecolor=red, linestyle=dashed, dash=1pt](2.6,0.2)(0.4,0.1){0}{180}\psellipticarc[linecolor=red](2.6,0.2)(0.4,0.1){-180}{0}

\rput(2.6,0){\small $A_3$}
\psellipticarc[linecolor=red, linestyle=dashed, dash=1pt](2.6,-.2)(0.4,0.1){0}{180}\psellipticarc[linecolor=red](2.6,-.2)(0.4,0.1){-180}{0}

\psellipticarc[linecolor=red, linestyle=dashed, dash=1pt](1,0.2)(1,0.1){0}{180}
\psellipticarc[linecolor=red](1,0.2)(1,0.1){180}{360}
\rput(1,0){\small $A_2$}
\psellipticarc[linecolor=red, linestyle=dashed, dash=1pt](1,-0.2)(1,0.1){0}{180}
\psellipticarc[linecolor=red](1,-0.2)(1,0.1){180}{360}

\end{pspicture}
\caption{Decomposition of $\Sigma_g$ into annuli $A_i$ and shrunken pants.}
\label{shrunkenPP_fig}
\end{figure}

Given $(n_1,\ldots,n_{3g-3})$, in each shrunken pair of pants $F$, draw the standard presentation of $\mathfrak{m}\cap F$ using the intersection numbers with $F$'s boundary curves as described previously.\\

For every $i=1,\ldots,3g-3$, the annulus $A_i$ has two boundary components $a_i^{+}$ and $a_i^{-}$ that are parallel copies of $a_i$ and are boundary curves of pants $F^{+}$ and $F^{-}$, respectively ($F^{+}$ and $F^{-}$ could be the same if $a_i$ is part of a folded pair of pants). Orient $a_i^{+}$ and $a_i^{-}$ as boundary curves (these orientations will be opposite).
Then, given $(t_1,\ldots,t_{3g-3})$, for every $i=1,\ldots,3g-3$: if $n_i\neq 0$, in the annulus $A_i$, connect the $k$-th chosen point on $\alpha_i^{-}$ (this is the half of $a_i^{-}$ belonging to the chosen hexagon) to the $(k+t_i)$-th chosen point of $\alpha_i^{+}$ by a straight edge\footnote{Think of a cylinder as the quotient of $\R\times [0,1]$ by $\Z$. Then, ``straight" means the image of a straight line in this universal cover.} in the cylinder. Here, the indices $k+t_i$ are considered mod $n_i$. {\it Note that since the orientations on $a_i^{+}$ and $a_i^{-}$ are opposite, this definition is symmetric and well-defined}. Figure~\ref{TwistP_fig} illustrates the outcome when $n_i=3$ and $t_i=2$.

\begin{figure}
\begin{pspicture}(-7,-1)(0,1.5)
\psset{unit=1cm}

\psellipticarc[ linestyle=dashed, dash=1pt](1,.5)(1.5,0.2){0}{180}
\psellipticarc(1,.5)(1.5,0.2){180}{360}
\psellipticarc[showpoints=true,arrowscale=2]{->}(1,.35)(1.5,0.2){-10}{0}

\psline(2.5,.5)(2.5,-.5)
\psline(-.5,.5)(-.5,-.5)

\psellipticarc[linestyle=dashed, dash=1pt](1,-0.5)(1.5,0.2){0}{180}
\psellipticarc(1,-0.5)(1.5,0.2){180}{360}
\psellipticarc[showpoints=true,arrowscale=2]{<-}(1,-.35)(1.5,0.2){180}{190}

\pscircle*(1,.3){.07}\pscircle*(1.5,.32){.07}\pscircle*(.5,.32){.07}
\pscircle*(1,-.7){.07}\pscircle*(1.5,-.68){.07}\pscircle*(.5,-.68){.07}

\rput(1,-1){\small $a^-$}\rput(1,.55){\small $a^+$}

\psline(.5,-.68)(1.5,.32)
\psline[linecolor=red](1,-.7)(2.5,0)\psellipticarc[linecolor=red, linestyle=dashed, dash=1pt](1,0)(1.5,0.2){0}{180}
\psline[linecolor=red](.5,.3)(-.5,0)

\psline[linecolor=blue](1.5,-.68)(2.5,-.2)
\psellipticarc[linecolor=blue, linestyle=dashed, dash=1pt](1,-.2)(1.5,0.2){0}{180}
\psline[linecolor=blue](1,.3)(-.5,-.2)

\end{pspicture}
\caption{The picture in an annulus when $n_i=3$ and $t_i=2$.}
\label{TwistP_fig}
\end{figure}

\begin{convention}\label{stdConvn=0}
If $n_i=0$, the standard convention is to assume $t_i\geq 0$ so that the multicurve $\mathfrak{m}$ includes $t_i$ parallel copies of $a_i$ in $A_i$.
\end{convention}

By connecting all the curve segments described above in each shrunken pair of pants and annulus, we obtain the multicurve $\mathfrak{m}$ associated with the DT coordinates $(n,t)$.

\medskip

While DT coordinates give an embedding of $\cS'_g$ into $\Z^{6g-6}$ and provide a straightforward recipe for drawing standard representatives of multicurves, the fact that the twist coordinates may take negative values makes them unsuitable for defining a degree function on the skein algebra. Moreover, these coordinates do not behave well under the product of multicurves in the skein algebra of the surface. We address the first issue by using geometric intersection numbers to coordinatize multicurves. We overcome the second issue by employing a modified version of Dehn--Thurston coordinates recently introduced in~\cite{BKL}; cf.~Section~\ref{Sec:proof-of-MainTheorem}. This approach requires an explicit understanding of the relationship between these coordinate systems, which is the subject of the next section.

\section{Thurston's compactification}\label{sec:TC}

In this section, we review the construction of Thurston's compactification and the $(9g-9)$-Theorem. Our main technical result, Theorem~\ref{f-functions}, which may be of independent interest to topologists and may have further applications in the study of skein algebras, provides an explicit formula for intersection numbers with a collection of $9g-9$ curves in terms of Dehn--Thurston coordinates. In the next section, we state a result concerning the behavior of these intersection numbers under the product map in the skein algebra.\medskip

In a nutshell, Thurston’s compactification of Teichm\"uller space adds a sphere at infinity of dimension $6g-7$, which serves as the base of a real cone identified with the space of measured foliations -- real-weighted analogues of multicurves on a surface. While this cone naturally lives in an infinite-dimensional space, it is often embedded into a finite-dimensional vector space, where it becomes a rational polyhedral cone over a piecewise-linear (PL) sphere. This PL sphere decomposes into a collection of convex polyhedra, each giving rise to a projective toric variety of twice the real dimension, and these varieties intersect along toric strata dictated by the structure of the polyhedral complex. Ultimately, we will identify the resulting toric complex with the boundary of a projective compactification of the character variety, thereby relating the sphere predicted by the geometric P=W conjecture to Thurston's boundary of Teichm\"uller space in a strong sense.\\

Suppose $\Sigma_g$ is a closed, oriented surface of genus $g$ with $g>1$. Recall that $\mathcal{S}_g$ denotes the set of isotopy classes of simple closed curves on $\Sigma_g$, and that $\mathcal{S}'_g$ denotes the set of isotopy classes of simple, closed, but not necessarily connected curves on $\Sigma_g$, called \textbf{multicurves}. For $\mf{m}_1, \mf{m}_2 \in \mathcal{S}'_g$, let $i(\mf{m}_1, \mf{m}_2)$ denote their geometric intersection number. 

\medskip

There are many ways to coordinatize multicurves and to use these coordinates to study the skein algebra.
In Section~\ref{sec:SK-DT}, we learned about Dehn-Thurston coordinates. Theorem~\ref{f-functions} below provides an explicit formula for transforming DT coordinates into certain intersection coordinates. There is also a modified version of DT coordinates that we will introduce and use later.\\

As we stated in the introduction, geometric intersection with elements of $\mathcal{S}_g$ induces a map
\bEq{iota_mc} 
\iota \colon \mathcal{S}'_g \lra \N^{\mathcal{S}_g} \subset \R_{\geq 0}^{\cS_g}
\eEq 
which is known to be injective.  It continuously extends to a similar embedding
\bEq{iota_F} 
\iota \colon \cF_g \lra  \R_{\geq 0}^{\cS_g}
\eEq
where $\mathcal{F}_g$ is the space of (Whitehead equivalence classes of) {\bf measured foliations} on $\Sigma_g$. There is a natural inclusion $\cS'_g\subset \cF_g$ so that the elements of $\cF_g$ can be thought of as multicurves with positive real weights.

\medskip

We will denote the projectivization of these embeddings, taking values in
$$
\mathbb{P}(\mathbb{R}_{\geq 0}^{\cS_g}) := \big( \mathbb{R}_{\geq 0}^{\cS_g} \setminus \{0\} \big) / \mathbb{R}_+,
$$
by $\iota_\mathbb{P}$. It is well known (cf.~\cite[Ch.~1]{FLP}) that:
\begin{itemize}
\item $\iota(\cF_g)$ is a real closed cone over $\iota_\mathbb{P}(\cF_g) \cong S^{6g-7} \subset \mathbb{P}(\mathbb{R}_{\geq 0}^{\cS_g})$.
\item Moreover,
\begin{equation} \label{closure}
\overline{\iota_\mathbb{P}(\cS_g)} = \overline{\iota_\mathbb{P}(\cS'_g)} = \iota_\mathbb{P}(\cF_g).
\end{equation}
\end{itemize}

The proof of the statements above in \cite{FLP} is by passing to finite dimension and considering a suitable finite set of curves that yields an embedding in the following sense.
\medskip

Fix a pants decomposition  
$$
\cP =\big\{a_1, \ldots, a_{3g-3}\big\}
$$  
of $\Si_g$ and an embedded dual graph $\Gamma$ as in Section~\ref{sec:SK-DT}. For each pants curve $a_j$, we consider a dual curve $a'_j$ (uniquely specified by the choice of $\Gamma$) that intersects $a_j$ but no other $a_i$, as illustrated in Figure~\ref{dual curve}. 
\begin{figure}
\begin{pspicture}(-12,-1.5)(0,.8)
\psset{unit=1cm}

\psellipse(-7,1)(0.4,0.1) 
\psellipse(-5,1)(0.4,0.1) 
\psarc(-6,1){.6}{180}{360}

\psellipse(-7,-1.5)(0.4,0.1) 
\psellipse(-5,-1.5)(0.4,0.1) 
\psarc(-6,-1.5){.6}{0}{180}

\psline(-7.4,1)(-7.4,-1.5)
\psline(-4.6,1)(-4.6,-1.5)

\psellipticarc[linestyle=dashed, dash=2pt](-6,-.25)(1.4,0.1){0}{180}
\psellipticarc(-6,-.25)(1.4,0.1){180}{360} \rput(-5.3,-.5){\small $a$}

\psellipticarc[linecolor=red, linestyle=dashed, dash=2pt](-6,-0.25)(0.2,.65){-90}{90} \rput(-6.3,.2){\small $a'$}
\psellipticarc[linecolor=red](-6,-0.25)(0.2,.65){90}{270}

\psarc(-1,0){1}{90}{270}
\pscircle[linecolor=red](-1.1,0){.5} \rput(-1.1,.67){\small $a'$}
\psline(-1,1)(0,1)
\psline(-1,-1)(0,-1)
\psarc(-1.5,0){.5}{-45}{45}
\psarc(-0.7,0){.5}{145}{215}
\psellipticarc[ linestyle=dashed, dash=1pt](-1.6,0)(0.4,0.1){0}{180} \rput(-1.8,-.25){\small $a$}
\psellipticarc(-1.6,0)(0.4,0.1){-180}{0} 
\psellipse(0,0)(0.2,1)

\end{pspicture}
\caption{Left: the dual curve $a'$ when $a$ belongs to two different pairs of pants. Right: the dual curve $a'$ when $a$ belongs to only one pair of pants.}
\label{dual curve}
\end{figure}
We let $a''_j$ to be the curve obtained from $+1$ Dehn twist\footnote{One may equally consider $-1$ Dehn twist or any other Dehn twist of $a_j'$ along $a_j$.} of $a'_j$ along $a_j$. Let $\cC$ denote the collection of $9g-9$ curves 
$$
\cC=\big(\cP, \cQ, \cQ_+\big) = \Big( (a_j)_{j=1}^{3g - 3}, (a'_j)_{j=1}^{3g - 3}, (a''_j)_{j=1}^{3g - 3} \Big).
$$

\begin{theorem}(\cite[Theorem~4.10]{FLP} or \cite[p.~438]{FM})\label{thm:ThurstonComp}
Taking intersections with the curves in $\cC$ (respectively, taking measures along the curves $\cC$) yields an embedding of $\cS'_g$ (respectively, 
$\cF_g$) into $\N^{9g-9}$ (respectively, $\R_{\geq 0}^{9g-9}$). Furthermore, the image of $\cF_g$ is a cone over $S^{6g-7}$.
\end{theorem}

As we explain in Remark~\ref{rmk-FLP} below, the proofs of~\cite[Theorem~4.10]{FLP} and~\cite[p.~438]{FM} are rather involved and do not establish an explicit relationship between Dehn--Thurston coordinates and the intersection coordinates associated with the collection $\cC$. \emph{Addressing this shortcoming, the following technical result provides explicit formulas for the intersection numbers with $a_j'$ and $a_j''$ in terms of Dehn--Thurston coordinates. We will need this explicit result to understand the image of $\cS'$ in $\N^{9g-9}$, to prove Theorem~\ref{ConeTheorem}, and ultimately to prove Theorem~\ref{MainTheorem}.}\\

 Given a pants decomposition and an embedded dual graph $\Gamma$, for every $a\in \cP$ let $a'\in \cQ$ and $a''\in \cQ_+$ be as above. For any multicurve $\mf{m}$ we are interested in formulas for $i(\mf{m},a')$ and $i(\mf{m},a'')$ in terms of DT coordinates of $\mf{m}$. We consider two cases: \medskip 
 
(I)  If $a \in \cP$ belongs to two pairs of pants, and  $a_1,a_2,a_3,a_4$ are the other four curves in $\cP$ that belong to the two pairs of pants containing $a$ as in Figure~\ref{dual curve2}-Left, let 
$$
\aligned
&\de_{14}=\max\big\{ \frac{n_1+n_4 -n}{2},0\big\},\quad\de_{23}=\max\big\{ \frac{n_2+n_3 -n}{2},0\big\},\\
&u_1=\max\big\{ \frac{n_1-n_4 -n}{2},0\big\},\quad u_4=\max\big\{ \frac{n_4-n_1 -n}{2},0\big\}\\
&u_2=\max\big\{ \frac{n_2-n_3 -n}{2},0\big\},\quad u_3=\max\big\{ \frac{n_3-n_2 -n}{2},0\big\}\\
&\wt{n}_1= n_1-2u_1-\de_{14},\quad \wt{n}_4= n_4-2u_4-\de_{14},\\
&\wt{n}_2= n_2-2u_2-\de_{23},\quad \wt{n}_3= n_3-2u_3-\de_{23},\\
\endaligned
$$
where $(n,n_1,n_2,n_3,n_4,n_5)$ are the geometric intersection numbers of $\mf{m}$ with $(a,a_1,a_2,a_3,a_4)$.\medskip

\begin{figure}
\begin{pspicture}(-12,-1.7)(0,.8)
\psset{unit=1cm}

\psarc(-6,1){.6}{180}{360}

\psarc(-6,-1.5){.6}{0}{180}

\psline(-7.4,1)(-7.4,-1.5)
\psline(-4.6,1)(-4.6,-1.5)

\psellipticarc[linestyle=dashed, dash=2pt](-6,-.25)(1.4,0.1){0}{180}
\psellipticarc(-6,-.25)(1.4,0.1){180}{360} \rput(-5.3,-.5){\small $a$}

\psellipticarc[linecolor=red, linestyle=dashed, dash=2pt](-6,-0.25)(0.2,.65){-90}{90} \rput(-6.3,.2){\small $a'$}
\psellipticarc[linecolor=red](-6,-0.25)(0.2,.65){90}{270} 

\psellipse(-7,1)(0.4,0.1) \rput(-7,1.3){\small $a_1$}
\psellipse(-5,1)(0.4,0.1) \rput(-5,1.3){\small $a_4$}
\psellipse(-7,-1.5)(0.4,0.1) \rput(-7,-1.8){\small $a_2$}
\psellipse(-5,-1.5)(0.4,0.1) \rput(-5,-1.8){\small $a_3$}

\psarc(-1,0){1}{90}{270}
\pscircle[linecolor=red](-1.1,0){.5} \rput(-1.1,.67){\small $a'$}
\psline(-1,1)(0,1)
\psline(-1,-1)(0,-1)
\psarc(-1.5,0){.5}{-45}{45}
\psarc(-0.7,0){.5}{145}{215}
\psellipticarc[ linestyle=dashed, dash=1pt](-1.6,0)(0.4,0.1){0}{180} \rput(-1.8,-.25){\small $a$}
\psellipticarc(-1.6,0)(0.4,0.1){-180}{0} 
\psellipse(0,0)(0.2,1)

 \rput(0,0){\small $a_1$}
 
\end{pspicture}
\caption{Left: case I. Right: case II}
\label{dual curve2}
\end{figure} 

(II) If $a\in\cP$ belongs to just one pair of pants, and  $a_1$ is the other curve in $\cP$ that belongs to the pairs of pants containing $a$ as in Figure~\ref{dual curve}-Right, define  
$$
u_1=\max\Big\{\frac{n_1}{2}-n,0\Big\},
$$
where $(n,n_1)$ are defined similarly.\medskip

The quantity $\delta_{14}$ is the number of model $\beta$-arcs in the standard presentation of $\mathfrak{m}$ joining the boundary components $a_1$ and $a_4$ (see the $\beta$-arcs in Figure~\ref{Triangular_fig}). Similarly, $\delta_{23}$ is the number of arcs joining the boundary components $a_2$ and $a_3$, and $u_i$ denotes the number of $U$-type arcs (see Figure~\ref{Dominant-xi_fig}) in the corresponding (shrunken) pair of pants with endpoints on the boundary component $a_i$. Finally, $\wt{n}_i$ is the number of model $\beta$-arcs joining the boundary component $a_i$ to $a^\pm$ in the corresponding (shrunken) pair of pants. Note that both $\wt{n}_1+\wt{n}_4$ and $\wt{n}_2+\wt{n}_3$ are less than or equal to $n$. The quantity $u_1$ in case~(II) has a similar meaning.

 \begin{theorem}\label{f-functions}
For every other multicurve $\mf{m}\in \cS'$ and $a\in \cP$, in case (I) above, the intersection numbers with $a'$ and $a''$ are given by piece-wise linear equations 
\bEq{FI}
n'=f(n,t,n_1,n_2,n_3,n_4) \quad \tn{and}\quad n''=f(n,t-n,n_1,n_2,n_3,n_4) ,
\eEq
where 
$$
\aligned
f(n,t,n_1,n_2,n_3,n_4)&=\max \Big\{ \abs{2t+n-\wt{n}_1-\wt{n}_3}, \wt{n}_1+\wt{n}_3-n, \wt{n}_2+\wt{n}_4-n\Big\}\\
&\,+\de_{14}+\de_{23}+2(u_1+u_2+u_3+u_4)
\endaligned
$$
\medskip
In this case (II), we have 
\bEq{FIfolded}
n'=h(n,t,n_1) \quad \tn{and}\quad n''=h(n,t-n,n_1) ,
\eEq
where $h(n,t,n_1)=|t|+u_1$
 \end{theorem}
 
While the formulas may appear intimidating, in Section~\ref{sec:DTandf-functions} we reduce the statements to simpler forms.

 \begin{remark}
 A similar result for the combinatorial Teichm\"uller space is proved in \cite{ABCGLW}; see Lemma~2.43 and 2.45. The formulas have similarities but are very different.
 \end{remark}
 
 \begin{remark}\label{rmk-FLP}
This remark explains some details of the proof of~\cite[Theorem~4.10]{FLP} (as well as~\cite[p.~438]{FM} for $\cF_g$) and, in doing so, clarifies the novelty of the result above.
\medskip

The intersection numbers directly used in the proofs of~\cite[Theorem~4.10]{FLP} (and~\cite[p.~438]{FM} for $\cF_g$) are not simply intersection numbers with closed curves in $\cC$. Instead, their arguments rely on an intermediate system of coordinates.
\medskip

Given $(\cP,\Gamma)$, $3g-3$ of these intermediate $9g-9$ coordinates are intersection numbers with the curves in $\cP$. The remaining coordinates are intersection numbers with certain local arcs in a thickening of the curves $a_j$.
\medskip

More precisely, the collection of curves $(K_j,K'_j,K''_j)_{j=1}^{3g-3}$ considered in~\cite[p.~63, (I)(II)]{FLP} coincides with the triple $(\cP,\cQ,\cQ_+)=(a_j,a'_j,a''_j)_{j=1}^{3g-3}$ used above. However, for an explicit closed cone $B$, the classifying maps
\[
\ov\Phi\colon \cF_g\lra B\subset \R^{9g-9}_{\geq 0}
\quad \text{and} \quad
\Phi\colon \cS'_g \lra B \cap \N^{9g-9}
\]
defined in~\cite{FLP} do not use measures along/intersection numbers with the curves $a_j'$ and $a''_j$ from $\cQ$ and $\cQ_+$. Instead, $\Phi$ is defined using intersection numbers with the curves $a_j$ in $\cP$, together with intersection numbers with local arcs $S_j$ and $T_j$ in thickenings of the curves $a_j$; see the definitions preceding~\cite[Lemma~4.7]{FLP}. The map $\ov\Phi$ is defined analogously. Intersection numbers with $S_j$ and $T_j$ satisfy simple relations in terms of the twist coordinates along $a_j$, making the map $\Phi$ readily expressible in terms of Dehn--Thurston coordinates.\\

In~\cite[Appendix~C]{FLP}, it is shown that the twisting parameters $t_i$ can be recovered from $(n'_i,n''_i)$, but no explicit formulas are given; see the proof of Proposition~6.12. Similarly, for the map $\ov\Phi$,~\cite[Appendix~C]{FLP} shows that measures along $(S_j,T_j)$ can be obtained from measures along $(a'_j,a''_j)$. This (cf.~\cite[Theorem~4.10]{FLP}) establishes the existence of a closed cone $C\subset \R_{\geq 0}^\cS$ and a continuous map
\[
\theta\colon C\to B
\]
that is positively homogeneous of degree one and makes the following diagram commute:
\[
\xymatrix{
\cF_g \ar[rd]_{\ov\Phi}\ar[rr]^\iota && C \ar[ld]^\theta \\
& B &
}
\]
Moreover, $\theta$ induces a homeomorphism between $\iota(\cF_g)\subset C$ and $B$. However, as noted in~\cite[Sec.~6.6]{FLP}, the restriction
\[
\ov\Phi|_{\cS'_g}\colon \cS'_g\lra B\subset \N^{9g-9}
\]
does not coincide with $\Phi$. This discrepancy arises because the definition of $\Phi$ involves intersection numbers with arcs rather than closed curves, leading to certain discontinuity issues. Consequently, one cannot use the (not very explicit) formulas for $\theta^{-1}$ from~\cite[Appendix~C]{FLP} to express intersection numbers with $(S_j,T_j)$ in terms of intersection numbers with $(a'_j,a''_j)$. The source of this difficulty is summarized in~\cite[p.~68]{FLP} as follows: ``\emph{It seems very difficult to make these last formulas explicit.}'' 
 \end{remark}

 We end this section with the proof of Theorem~\ref{ConeTheorem} which is a refinement of Theorem~\ref{thm:ThurstonComp}.\\

 {\bf Proof of Theorem~\ref{ConeTheorem}.}
By (\ref{closure}) (see also \cite[Prp.~5.11]{FLP}), the formulas of Theorem~\ref{f-functions} extend to identical formulas for measured foliations.
As we stated in Theorem~\ref{thm:ThurstonComp} and explained in the remark above, the result that taking intersection numbers with the curves in $\mc{C}=(\cP,\cQ,\cQ_+)$ yields embeddings 
$$
\aligned
& \iota_\cC \colon \mathcal{S}'_g \lra \N^{9g-9}\\
&\iota_\cC \colon \mathcal{F}_g \lra \R_{\geq 0}^{9g-9}.
\endaligned
$$
is proved in \cite{FLP} via comparison with another embedding that comes from taking intersections with the curves in $\cP$ and some local arcs in the cylinders $A_i$. The formulas of Theorem~\ref{f-functions}, however, give a direct proof of the embedding property. They also give a more concrete description of how the image $\La=\iota_\cC\big(\mathcal{F}_g\big)$ looks like.\\

To prove that $\iota_\cC$ is one-to-one, we need to show that the twist parameter $t$ along a pants curve $a$ can be recovered from $(n, n', n'', n_1, \ldots, n_4)$ in case (I), and from $(n, n', n'', n_1)$ in case (II) of Theorem~\ref{f-functions}. In other words, fixing $(n, n_1, \ldots, n_4)$, we want to show that the map $t \mapsto (n'(t), n''(t))$ is an embedding $\R \lra \R_{\geq 0}^2$. We will demonstrate this both geometrically and algebraically.
 \\
 
In case (I), the formulas for $n'(t)$ and $n''(t)$ are of the form 
 $$
 n'(t)=\tn{max}\lrc{ |2t-A|,A,B}+C\quad \tn{and}\quad n''(t)=\tn{max}\lrc{ |2(t-n)-A|,A,B}+C
 $$
 where $A$, $B$, and $C$ are functions of $(n, n_1, \ldots, n_4)$. It is also easy to see that $A,B\leq n$ and $A+B\leq 0$. \medskip
 
 In case (II), 
  the formulas for $n'(t)$ and $n''(t)$ are of the form 
 $$
 n'(t)= |t|+C\quad \tn{and}\quad n''(t)=|t-n|+C
 $$
 where  $C$ is a function of $(n, n_1)$; thus, it is special case of case (I).\\
 
 The graphs of $n'$ and $n''$ (in their most general form) look as in Figure~\ref{graph-of-n'}; i.e. the interiors  of the intervals on which $n'$ and $n''$ are the constant $\tn{max}\{A,B\}+C$ have no overlap and the image of the map $t\mapsto(n'(t),n''(t))$ (in its most general form) is an embedding as in Figure~\ref{nnembedding}. \\

\begin{figure}
\begin{pspicture}(-7,.5)(0,3)
\psset{unit=1.2cm}

\psline(-1,1)(-2,3)
\psline(1,1)(2,3)
\psline(-1,1)(1,1)
\psline{<->}(-1,.7)(1,.7)\rput(0,.5){\small{at most $n$}}
\psline{->}(0,.8)(2.5,.8)\rput(2,.6){\small{shifted by $n$}}
\rput(0,1.5){\small{$n'(t)$}}
\rput(2.5,1.5){\small{$n''(t)$}}

\psline(1.5,1)(.5,3)
\psline(3.5,1)(4.5,3)
\psline(1.5,1)(3.5,1)

\end{pspicture}
\caption{Graph of $n'$.}
\label{graph-of-n'}
\end{figure}

\begin{figure}
\begin{pspicture}(-5,1)(0,5)
\psset{unit=1cm}

\psline{<->}(5,3)(3,1)(2,1)(1,2)(1,3)(2.5,4.5)
\rput(5.5,3.3){\small{$t\to \infty$}}
\rput(3.3,4.5){\small{$t\to -\infty$}}

\end{pspicture}
\caption{Image of the parametric curve $(n'(t),n''(t))$}
\label{nnembedding}
\end{figure}

 Algebraically, in case (I),  twisting parameter $t$ is the following function of the intersection numbers:
 \begin{equation}\label{eq:t-equation}
 t=\begin{cases}
n- (n''-A-C)/2 &\tn{if}\quad  n'\leq n''\\
(n'+A-C)/2& \tn{if}\quad  n'\geq n''\;.
 \end{cases}
 \end{equation}
 In case (II), we get the simpler equations
$$
 t=\begin{cases}
n- (n''-C) &\tn{if}\quad  n'\leq n''\\
n'-C& \tn{if}\quad  n'\geq n''\;.
\end{cases}
$$

In the special case of $n=0$, the formulas for $n'(t)$ and $n''(t)$ in cases (I) and (II) simplify to 
 $$
 \aligned
 &n'(t)=n''(t)=|2t|+C,\\
 &n'(t)=n''(t)= |t|+C,
 \endaligned
 $$
 respectively.  Therefore, the standard convention (i.e. Convention~\ref{stdConvn=0}) requires taking $t=\frac{n'-C}{2}\geq 0$ in case (I) and $t=n'-C\geq 0$ in case II as the twist coordinate; see more about this in Remark~\ref{Rmkont} below. \\

In conclusion, having parametrized the multicurves with DT coordinates, the embedding  $\iota_\cC$ is the piecewise-linear embedding
\begin{equation}\label{iotaC}
\aligned
&\mc{D}_g=\lrc{(n,t)\in \R^2\colon n\geq 0, ~~t\geq 0~\tn{if}~n=0}^{3g-3} \lra \R_{\geq 0}^{9g-9},\\
&(n_1,\ldots,n_{3g-3},t_1,\ldots,t_{3g-3})\lra \big(n_1,\ldots,n_{3g-3}, n'_1(n,t), \ldots , n'_{3g-3}(n,t),n''_1(n,t), \ldots , n''_{3g-3}(n,t)\big)
\endaligned
\end{equation}
defined by explicit equations in Theorem~\ref{f-functions}.\medskip

The cone 
$$
\La=\iota_\cC(\cF_g)\subset \R_{\geq 0}^{9g-9}
$$ is a closed cone, while the domain $\mc{D}_g$ of $\iota_\cC$ is a half-open subcone of $\R^{6g-6}$. One can also turn $\mc{D}_g$ into a closed singular flat manifold
$$
(\R^2/\pm)^{3g-3}
$$ 
by identifying the (included in $\mc{D}_g$) positive $t_i$-axis (i.e. $n_i=0$ and $t_i\geq 0$) and the (excluded in $\mc{D}_g$) negative $t_i$-axis in $\R_{\geq 0}\times \R$ of the closed cone 
$$
\ov{\mc{D}_g}=(\R_{\geq 0} \times \R)^{3g-3}.
$$
The map \eqref{iotaC} continuously extends to a piecewise linear map
$$
\ov{\mc{D}_g}\lra  \R_{\geq 0}^{9g-9}
$$
that identifies the equivalent points in $\partial \ov{\mc{D}_g}$; i.e. $\iota_\cC$ can be seen as a piecewise linear embedding
$$
(\R^2/\pm)^{3g-3}\lra \R_{\geq 0}^{9g-9}
$$
between the two closed cones.\medskip

Once again, $\iota_\cC$ is a homogeneous piecewise-linear function because $(n',n'')$ are defined using absolute values and maxima of collections of linear functions. The closed domain $\ov{\cD}_g$ decomposes as a union of closed rational polyhedral cones $\La'_\ell$, on each of which $\iota_\cC$ is linear. These regions are obtained by determining the different branches of linearity of the maximum and absolute value functions, which result in inequalities on the DT--coordinates $(n,t)$ specifying each $\La'_\ell$. More specifically, we restrict to domains on which all formulas for $\de_{ij}$, $u_i$, $|2t-A|$, $\tn{max}\lrc{|2t-A|,A,B}$, and $\tn{max}\lrc{|2(t-n)-A|,A,B}$ become linear. We will refer to these inequalities abstractly in the proofs as
\[
I_{\ell,1}\geq 0,\ldots,I_{\ell,k_\ell}\geq 0,
\]
and they are described explicitly below.\medskip

The image $\La_\ell$ of each $\La'_\ell$ under $\iota_\cC$ is then a rational polyhedral cone contained in $\La=\iota_\cC(\cF)$. The cones $\La_\ell$ decompose $\La$ into a union of finitely many rational polyhedral cones glued along their faces.\medskip

The inequalities mentioned above are
\[
n_i\geq n_j+n_k \quad\tn{or}\quad n_i\leq n_j+n_k
\]
in each pair of pants $F(a_i,a_j,a_k)$, and
\[
A=\wt{n}_1+\wt{n}_3-n \geq 0 \ \tn{or}\ \leq 0
\quad\tn{and}\quad
B=\wt{n}_2+\wt{n}_4-n \geq 0 \ \tn{or}\ \leq 0
\]
in every adjacent pair of shrunken pairs of pants, as in Theorem~\ref{f-functions}. The first set of inequalities determines the values of the max functions defining $\de_{ij}$ and $u_i$, and hence $\wt{n}_i$, $A$, $B$, and $C$. The second set, together with inequalities on $t$ described below, is required to determine whether $\tn{max}\lrc{|2t-A|,A,B}$ can be equal to constants (in $t$) $A$ or $B$. Once the inequalities on intersection numbers are fixed, the cones $\La'_\ell$ are determined by choosing intervals for the twist parameters $t$ that make $n'$ and $n''$ linear in the twist (and intersection) coordinates; see Section~\ref{sec:g2} for some genus--two examples. The following cases exhaust all possibilities:
\begin{enumerate}
\item If both $A=\wt{n}_1+\wt{n}_3-n$ and $B=\wt{n}_2+\wt{n}_4-n$ are non-positive, then $n'(t)$ and $n''(t)$ simplify to
\[
n'(t)=|2t-A|+C \quad\tn{and}\quad n''(t)=|2(t-n)-A|+C.
\]
In this case, there are three domains of linearity for $t$:
\[
t\leq \frac{A}{2},\qquad
\frac{A}{2}\leq t \leq n+\frac{A}{2},\qquad
\tn{and}\quad t\geq n+\frac{A}{2}.
\]
On these three regions, $n'$ and $n''$ satisfy the linear relations $n''-n'=2n$, $n'-n''=2n$, and $n'+n''=2n+2C$, respectively.
\item If exactly one of $A$ or $B$ is positive (which must be the case, since $A+B<0$), then there are five domains of linearity for $t$.
\begin{enumerate}
\item If $A>0>B>-A$, these five intervals are
\[
t\leq 0,\qquad A\leq t \leq n,\qquad t\geq n+A,
\]
and
\[
0\leq t\leq A,\qquad n\leq t\leq n+A.
\]
\item If $-B<A<0<B$, they are
\[
t\leq \frac{A-B}{2},\qquad
\frac{A+B}{2}\leq t \leq n+\frac{A-B}{2},\qquad
t\geq n+\frac{A+B}{2},
\]
and
\[
\frac{A-B}{2}\leq t\leq \frac{A+B}{2},\qquad
n+\frac{A-B}{2}\leq t\leq n+\frac{A+B}{2};
\]
see Figure~\ref{graph-of-n'2}.
\end{enumerate}
\begin{figure} \begin{pspicture}(-7,0.2)(0,3) \psset{unit=1.2cm} \psline(-1,1)(-2,3) \psline(1,1)(2,3) \psline(-1,1)(1,1) \psline(-1,.6)(-1,.8)\rput(-1,.3){\tiny{$\frac{A-B}{2}$}} \psline(1,.6)(1,.8)\rput(.7,.3){\tiny{$\frac{A+B}{2}$}} \psline(1.5,.6)(1.5,.8)\rput(1.7,.3){\tiny{$n+\frac{A-B}{2}$}} \psline(3.5,.6)(3.5,.8)\rput(3.7,.3){\tiny{$n+\frac{A+B}{2}$}} \psline{<->}(-3,.7)(6,.7)\rput(-3,.3){\tiny{$-\infty$}}\rput(6,.3){\tiny{$\infty$}} \rput(0,1.5){\small{$n'(t)$}} \rput(2.5,1.5){\small{$n''(t)$}} \psline(1.5,1)(.5,3) \psline(3.5,1)(4.5,3) \psline(1.5,1)(3.5,1) \end{pspicture} \caption{Domains of linearity of $(n'(t),n''(t))$ when $-B<A<0<B$.} \label{graph-of-n'2} \end{figure}

In each subcase of~(2), on the first three intervals the relations between $n'$ and $n''$ are as before. On the remaining two intervals, one of $n'$ or $n''$ is constant as a function of $t$, while the other is linear in $2t$ or $-2t$.\medskip
\end{enumerate}
In the case $g=2$, only Case~(1) occurs, yielding three intervals for the twist parameter; see Section~\ref{sec:g2}.\medskip

Summarizing, each cone $\La'_\ell$ is a convex rational polyhedral cone in $\ov{\mc{D}_g}$ defined by a collection of inequalities, and $\La_\ell$ is the image of $\La'_\ell$ under a linear embedding
\[
L_\ell \colon \R^{6g-6}\lra \R^{9g-9}.
\]
The cone $\La_\ell$ itself can be described by a system of $\Q$--linear inequalities and equalities inside $\R_{\geq 0}^{9g-9}$, although, as suggested by the relations between $n'$ and $n''$ above, these descriptions are somewhat cumbersome to write explicitly.\medskip

By Thurston's result, $\La$ is a cone over $S^{6g-7}$. Thus, each (maximal) $\La_\ell$ is a cone over one of the $(6g-7)$--dimensional polytopes $P_\ell$ in a polyhedral decomposition 
$$
P=\bigcup_\ell P_\ell =  (\La-\{0\})/\R_{+}
$$ 
of a piecewise-linearly embedded
\[
S^{6g-7}\subset \De_{9g-10}
=\big\{(x_1,\ldots,x_{9g-9})\colon x_i\geq 0,\ \sum x_i=1\big\}.
\]
It follows from the formulas in Theorem~\ref{f-functions}, \eqref{Pants-coordinates_e}, and \eqref{eq:t-equation} that $\iota_\cC(\cS'_g)$ contains $\La\cap 4\N^{9g-9}$. Consequently, the image of $\cS'_g$ in each $\La_\ell$ is a finite-index submonoid of $\La_\ell\cap \N^{9g-9}$.

\begin{remark}\label{Rmkont}
For each pants curve, the closed domain $\ov{\cD}_g$ contains the entire real line $\R$ as the range of possible values of $t$. Which half-line, $\R_{\geq 0}$ or $\R_{\leq 0}$, belongs to $\La'_\ell$ depends on whether $n'\geq n''$ or $n''\geq n'$, in the following sense.\medskip

A multicurve satisfying $n=0$ lies in the intersection of two distinct cones $\La_\ell$, locally distinguished by the inequalities $n'\geq n''$ and $n''\geq n'$, which intersect along $n'=n''$. By \eqref{eq:t-equation}, approaching $n'=n''$ from the side where $n'\geq n''$ yields the non-negative limit
$$
t=\frac{n'-C}{2}.
$$
Approaching from the side where $n''\geq n'$ yields the non-positive limit
$$
t=-\frac{n'-C}{2}.
$$

These two possibilities are identified under the antipodal map 
$$
(\R^2/\pm)^{3g-3}=\ov{\mc{D}}_g/\sim,\quad
(n,t)\sim (n',t') \Longleftrightarrow \forall\, i=1,\ldots,3g-3 :
\begin{cases}
(n_i,t_i)=(n'_i,t'_i), & \text{if } n_i\neq 0,\\
t_i=\pm t'_i, & \text{if } n_i=n'_i=0,
\end{cases}
$$
that glues the two cones along a common face.
\end{remark}

\begin{convention}\label{new-convention}
The discussion in the remark above shows that, for each $\La_\ell$ and each pants curve $a\in\cP$, depending on whether $n'\geq n''$ or $n''\geq n'$ on $\La_\ell$ near a boundary point with $n=0$, Convention~\ref{stdConvn=0} for the twist coordinate on $\La'_\ell$ along $a$ when $n=0$ must be modified to require either $t\geq 0$ or $t\leq 0$; respectively. In other words, the convention for a multicurve $\mf{m}$ should depend on the cone in which $\mf{m}$ is regarded as an element.
\end{convention}

We summarize the discussion above in the following items:
\begin{itemize}
\item Thinking of DT coordinates as points in 
$$
(\R^2/\pm)^{3g-3},
$$ 
there is a decomposition 
$$
(\R^2/\pm)^{3g-3}=\bigcup_{\ell} \La'_\ell,
$$
where each $\La'_\ell$ can be identified with a convex rational polyhedral cone in $\ov{\mc{D}}_g=(\R_{\geq 0} \times \R)^{3g-3}$, and $\La_\ell$ is the image of $\La'_\ell$ under a $\Q$-linear embedding
$$
L_\ell=\iota_\cC|_{\La'_\ell}\colon \La'_\ell\lra \La_\ell\subset \R_{\geq 0}^{9g-9}.
$$ 
\item Each $\La'_\ell$ is determined by, (1) a set of inequalities on intersection numbers $n$, and (2) by a set of inequalities on twist parameters $t$ with bounds depending on $n$.
\item For a multicurve $\mf{m}$ having image in the cone $\La_\ell$, if the intersection number $n$ of $\mf{m}$ with a pants curve $a$ is $n=0$, we modify Convention~\ref{stdConvn=0} so that its twist number $t$ takes values in $\R_{\geq 0}$ if $n'\geq n''$, and in $\R_{\leq 0}$ if $n'\leq n''$.
\end{itemize}
\medskip

In order to show that the image $Q_\ell$ of $\cS'_g$ in each $\Lambda_\ell$ is a saturated monoid (cf. Definition~\ref{monoids_dfn}), it is more convenient -- as is known in the punctured case (cf.~\cite[Proposition~3.2]{AF}) -- to work with the so-called {\bf angle} coordinates, defined as follows. We will also use these angle coordinates in the computations of Section~\ref{sec:g2}.
\medskip

In each pair of pants $F=F(a_1,a_2,a_3)$, as shown in Figure~\ref{Dominant-xi_fig}, define the angle coordinates $(\theta_1,\theta_2,\theta_3)$ by
$$
 \theta_{F,i} = \frac{n_j+n_k-n_i}{2}\in \Z\qquad \forall~i=1,2,3,\ \{j,k\}=\{1,2,3\}\setminus\{i\}.
$$
Note that for folded pairs of pants $F=F(a_1,2a)$, as in Theorem~\ref{f-functions}.(II), the (only) angle coordinate $\theta$ corresponding to $a$ satisfies
\bEq{theta-equality}
\theta=\frac{n_1}{2}\in \N.
\eEq
The total number of angle and twist coordinates is at most $9g-9$. By \eqref{Pants-coordinates_e} and Convention~\ref{new-convention}, the angle coordinates are integers subject to the following relations:
\begin{itemize}
\item With notation as in Theorem~\ref{f-functions}.(I), if
$$
(\theta_{F^+},\theta_{F^+,1},\theta_{F^+,4})\quad \text{and}\quad (\theta_{F^-},\theta_{F^-,2},\theta_{F^-,3})
$$
are the angle coordinates corresponding to $(a,a_1,\ldots,a_4)$ in adjacent pairs of pants $F^+=F(a,a_1,a_4)$ and $F^-=F(a,a_2,a_3)$, then
\bEq{theta-equality}
\theta_{F^+,1}+\theta_{F^+,4}=\theta_{F^-,2}+\theta_{F^-,3}.
\eEq
\item In each pair of pants $F=F(a_1,a_2,a_3)$,
\bEq{positivity-inequality}
n_i=\theta_{F,j}+\theta_{F,k}\geq 0\qquad \forall~i=1,2,3,\ \{j,k\}=\{1,2,3\}\setminus\{i\}.
\eEq
\item In each pair of pants $F=F(a_1,a_2,a_3)$, if $n_i=\theta_{F,j}+\theta_{F,k}=0$, then, depending on the chosen convention in Convention~\ref{new-convention}, we have $t_i\geq 0$ or $t_i\leq 0$.
\end{itemize}
Since all these conditions are equalities or inequalities, the integer points in each $\Lambda'_\ell$ corresponding to multicurves form a saturated monoid.
More precisely, the integer points in each $\Lambda'_\ell$ are the $\Z$-linear images of all the integer points in a rational polyhedral cone $\La''_\ell\subset \Z^N$, where $N$ is the total number of angle and twist coordinates, and $\La''_\ell$ is the cone specified by the conditions above together with the additional inequalities
$I_{\ell,1}\geq 0,\ldots,I_{\ell,k_\ell}\geq 0$
used to specify $\La'_\ell$
\medskip

Similarly, it is easy to see that the equations of Theorem~\ref{f-functions}.(II) are also $\Z$-linear in the angle and twist coordinates. Therefore, the image of $\cS'_g$ in each $\Lambda_\ell$ is the $\Z$-linear image of $\La''_\ell\cap \Z^N$. Since the latter is saturated, so is the former.\medskip

The following diagram illustrates the relation between these coordinate systems.

\begin{pspicture}(-5,-1)(9,2.4)
  \psellipse(-1,0)(1.4,0.8)
  \rput(-1,.2){\scriptsize{angle/twist}}
  \rput(-1,-.2){\scriptsize{coordinates} $ (\theta,t)$}
  
  \psellipse(4,0)(1.4,0.8)
    \rput(4,.2){\scriptsize{DT}}
  \rput(4,-.2){\scriptsize{coordinates} $(n,t)$}

  \psellipse(9,0)(1.4,0.8)
    \rput(9,.3){\scriptsize{intersection}}
  \rput(9,-.05){\scriptsize{coordinates}}
    \rput(9,-.4){\scriptsize{$(n,n',n'')$}}

  \psline{->}(0.5,0)(2.5,0)     \rput(1.5,-.2){\scriptsize{$\Z$-linear}}\rput(1.5,.2){\scriptsize{Piece-wise}}
  \psline{->}(5.5,0)(7.5,0)    \rput(6.3,-.2){\scriptsize{$\Q$-linear}}\rput(6.3,.2){\scriptsize{Piece-wise}}
  
    \pnode(0.2,.5){Aout}
  \pnode(2.8,0){Bin}
  \pnode(5.4,0){Bout}
  \pnode(7.8,.5){Cin}
  \nccurve[arrows=->,angleA=60,angleB=120,ncurv=.6]{Aout}{Cin}
  \rput(4,1.6){\scriptsize{Piece-wise $\Z$-linear}}
\end{pspicture}
\medskip

This finishes the proof of Theorem~\ref{ConeTheorem} and establishes several related facts that will be used in the subsequent sections.

\qed

\begin{remark}\label{hf-in-thea}
With notation as in Theorem~\ref{f-functions}, let $F^\pm$ denote the upper and lower pairs of pants. Then, the functions $f$ and $h$ take the following form with respect to angle/twist coordinates:
$$
f=\max \Big\{ \abs{2t-A}, A, B\Big\}\ + C
$$
where 
$$
\aligned
&\de_{14}=\max\big\{\theta_{F^+,a},0\big\},\quad\de_{23}=\max\big\{ \theta_{F^-,a},0\big\},\\
&u_1=\max\big\{-\theta_{F^+,a_1},0\big\},\quad u_4=\max\big\{-\theta_{F^+,a_4},0\big\}\\
&u_2=\max\big\{-\theta_{F^-,a_2},0\big\},\quad u_3=\max\big\{-\theta_{F^-,a_3},0\big\}\\
&A= (n_1+n_3-n)-2(u_1+u_3)- (\de_{14}+\de_{23}),\\
&B= (n_2+n_4-n)-2(u_2+u_4)- (\de_{14}+\de_{23}),\\
&C=(\de_{14}+\de_{23})+2(u_1+u_2+u_3+u_4)
\endaligned
$$
and 
$$
h=|t|+\max\Big\{-\theta_{F,a_1},0\Big\}.
$$

\end{remark}

\section{Compactification of $\cX_g$}\label{sec:conenction}
In this section, we use
\begin{itemize}
\item the isomorphism between the ring of regular functions $\C[\cX_g]$ and the classical skein algebra $\tn{Sk}_g$, and
\item the filtration on the skein algebra arising from intersection numbers with the curves of the collection $\cC$ in Theorem~\ref{ConeTheorem}
\end{itemize}
to construct a projective compactification of $\cX_g$ with toric boundary divisors.\\

Fix a collection $\cC$ as in Theorem~\ref{ConeTheorem}. The sum of intersection numbers with the curves in $\cC$ defines a degree function
\begin{equation}\label{grading}
|\mf{m}|\defeq \sum_{c\in \cC}  i(\mf{m},c)
\end{equation} 
on the multicurves in $\cS'_g$ which are the vector space generators of skein algebra $\tn{Sk}_g$.
This degree function in turn defines a filtration on the skein algebra
$$
 F_0\subset F_1 \subset \cdots \subset \tn{Sk}_g=\bigcup_{d\geq 0} F_d
$$
such that $F_d$ is the vector space generated by the set of multicurves $\mf{m}$ satisfying $|\mf{m}|\leq  d$.
To every such filtration, we can associate two graded algebras: 
\bEq{GAlg_e}
\tn{Sk}_{g}[x]\defeq \bigoplus_{d\geq 0} F_d x^d \quad \tn{and}\quad \tn{Sk}_{g}^{\tn{gd}}\defeq \bigoplus_{d\geq 0} F_d/F_{d-1}.
\eEq
where $x$ is an indeterminate over $\tn{Sk}_{g}$. We will denote the induced multiplication in the associated graded algebra $\tn{Sk}_{g}^{\tn{gd}}$ by $\overline{\ast}$. \\

By \cite[Proposition~2.8]{Mon} and the isomorphism $\tn{Sk}_g\cong \C[\cX_g]$, for every such filtration
\begin{equation}\label{eq:partialX}
\overline{\cX_g} \defeq \tn{Proj}(\tn{Sk}_{g}[u])
\end{equation}
is a reduced scheme containing $\cX_g$ as a Zariski dense open subset. The complement 
$$
D_g\defeq\partial  \ov\cX_g=\ov\cX_g -\cX_g
$$ 
is the zero set of the ideal $\ll x\rr$ generated by $x$ and is isomorphic (as a scheme) to 
$$
\tn{Proj}(\tn{Sk}_{g}^{\tn{gd}}).
$$ 
\begin{definition}
In general, we call a filtration on a $\C$-algebra  $\cR$  \textbf{finitely-generated} if the graded algebra $\cR[u]$ is finitely generated. We call it \textbf{projective} if it is finitely generated and $F_0=\C$.
\end{definition}
 If the filtration is finitely generated, then $\ov\cX_g $ is a closed sub-scheme of a {\it generalized weighted projective space}\footnote{Following \cite[Definition~2.7]{Mon}, a generalized weighted projective space is a weighted projective space possibly with negative weights.} and $\partial \ov\cX_g $ is the support of an effective ample divisor. If in addition the filtration is projective, then $\ov\cX_g $ is a projective variety; see \cite{Mon}.  The filtration coming from taking intersection with the collection of curves in Theorem~\ref{ConeTheorem} has these properties.    \\

In what follows we give a clear description of $\tn{Sk}_{g}^{\tn{gd}}$. For every two multicurves $\mf{m}$ and $\mf{m}'$, recall that 
\bEq{decomposition}
\mf{m}\ast\mf{m}'=\sum_{s} \la_s \mf{m}_s,
\eEq
where $s$ runs over all possible smoothings of $\mf{m}\cup \mf{m}'$ (in general position) and the coefficient  $\la_s$ collects the contributions of trivial components (trivial curve $=-2$) appearing in a smoothing. Passing to $\tn{Sk}_{g}^{\tn{gd}}$, we get 
\bEq{eq:ovast}
\mf{m}~\ov\ast~\mf{m}'=\sum_{s\colon |\mf{m}_s|=|\mf{m}|+|\mf{m}'|} \la_s \mf{m}_s,
\eEq

\medskip

The following theorem shows that the degree function defined with respect to the collection of $9g-9$ curves in Theorem~\ref{ConeTheorem} behaves well under the product  $\ast$ in the skein algebra. More specifically, it shows that the embedding of Theorem~\ref{ConeTheorem}  induces an isomorphism between the graded coordinate ring of $\cX_{g}$ and the Stanley–Reisner algebra $\C[Q]$ associated to  
$$
Q=\iota(\mathcal{S}'_g) \subset \Lambda.
$$ 
From Theorem~\ref{ThmC}  we will conclude that $\partial \ov\cX_g$ is a union of toric varieties intersecting along toric strata according to the combinatorics of the polytope complex $P$.

\begin{theorem}\label{ThmC} 
With notation as in Theorem~\ref{ConeTheorem} and \eqref{eq:ovast}, and for every two multicurves $\mf{m}$ and $\mf{m}'$, the followings hold.
\begin{enumerate}
\item If $\iota_\cC(\mf{m})$ and $\iota_\cC(\mf{m}')$  belong to the same cone $\Lambda_\ell$, for some $\ell$, then there exists a unique term on the righthand side of (\ref{decomposition}), say $\mf{m}_\#$, satisfying $ \iota_\cC(\mf{m}_\#)\in \Lambda_\ell$ and
$$
\mf{m}\,\overline{\ast}\,\mf{m}'= \mf{m}_\#.
$$
\item If $\iota_\cC(\mf{m})$ and $\iota_\cC(\mf{m}')$ do not belong to the same cone $\Lambda_\ell$, for any $\ell$, then 
$$
\mf{m}\,\overline{\ast}\,\mf{m}'=0\in \tn{Sk}_{g}^{\tn{gd}}.
$$ 
\end{enumerate}
\end{theorem}

The key idea behind the proof is as follows. As mentioned earlier, Dehn–Thurston coordinates do not behave well under the product map ${\ast}$ in $\mathrm{Sk}_{g}$. In \cite{BKL}, the authors introduced a modified version of DT coordinates that satisfy the following “leading order term property” (see \cite[Theorem~9.1]{FKL23}):

\begin{quote}
For any two multicurves $\mf{m}$ and $\mf{m}'$, there is a unique term $\mf{m}_\#$ in their product whose modified DT coordinates are equal to the sum of the modified DT coordinates of $\mf{m}$ and~$\mf{m}'$.
\end{quote}

Modified DT coordinates are piecewise linear functions of the standard DT coordinates. Our key observation is that, after minor adjustments, the modified DT coordinates and the intersection coordinates on each cone $\Lambda_\ell$ are related by linear equations. This yields the first statement. The second follows from a convexity argument. We postpone the full proof to Section~\ref{Sec:proof-of-MainTheorem}.

\begin{remark}
There is a different method to endow the space $\cF_g$ with the structure of a finite union of rational polyhedral cones glued along common faces, by considering (maximal birecurrent) train tracks and the closed cones of multicurves supported on such train tracks; see \cite[Section~2.6]{PH}. However, even in the genus--two case, it seems that the resulting cones are mostly different from those in Theorem~\ref{ConeTheorem} and do not enjoy the crucial second property of Theorem~\ref{ThmC}.
\end{remark}

{\bf Proof of Theorem~\ref{MainTheorem}.} 
First, we make a short digression and state a few basic facts about monoids and toric varieties.

\begin{definition}\label{monoids_dfn}
A monoid $M$ is called 
\begin{itemize}
\item {\bf integral}: if the map $M\lra M^{gp}$ is injective, where $M^{gp}$ is the group associated to $M$;
\item {\bf fine}: if $M$ is finitely generated and integral;
\item {\bf saturated}: if it is integral and for any $m\in \N$ and $q\in M^{gp}$ such that $mq\in M$, we also have $q\in M$;
\item {\bf sharp}: if $M$ has no invertible elements other than $0$;
\item {\bf toric}: if $M$ is fine-saturated and $M^{gp}$ is torsion-free (and thus, free of finite rank).
\end{itemize}
\end{definition}
It is fairly simple to show that if $M$ is fine, saturated, and sharp then it is toric.\medskip

In general, suppose  $\La \subset \R_{\geq 0}^k$ is a rational polyhedral fan; that is, $\La$ is a (finite) union $\La=\bigcup_\ell \La_\ell $ of convex rational polyhedral cones glued along common faces. In our case of interest, $k$ will be $9g-9$ and $\La$ is the fan consisting of the finite union of the cones in Theorem~\ref{ConeTheorem}. In any case, let $\La(\N)$ denote the subset $\La \cap \N^k$ of integer points.
Then, the Stanley–Reisner ring $\C[\La(\N)]$ is the free $\C$-vector space over the monomials $x^{\mf{n}}$, $\mf{n}\in \La(\N)$, where the product structure is defined by 
$$
x^{\mf{n}}\cdot x^{\mf{n}'}
=
\begin{cases}
x^{\mf{n}+\mf{n}'} & \tn{if there is a cone in } \La \tn{ that contains } \mf{n},\mf{n}'\\
0                              & \tn{otherwise}.
\end{cases}
$$
This is a graded ring where each monomial $x^{\mf{n}}$ is graded by the sum of coordinates of $\mf{n}$.
For each cone $\La_\ell \subset \La$, $\La_\ell(\N)$ is fine, saturated, and sharp; therefore, it is toric. Thus, we obtain a projective toric  variety 
$$
X_\ell=\tn{Proj}(\C[\La_\ell(\N)])
$$ 
with the moment polytope $P_\ell$ such that $\La_\ell$ is a cone over $P_\ell$. More precisely, $P_\ell$ is the intersection of $\La_\ell$ and the affine hyperplane $|\mf{n}|=1$. Furthermore, the projective variety
\bEq{Toric-space}
X=\tn{Proj}(\C[\La(\N)]) = \bigcup_{\ell} X_\ell
\eEq
associated to the entire $\La(\N)$ is a stratified space with locally closed strata indexed by subcones of $\La$; it is a union of toric varieties meeting along toric strata according to combinatorics of the moment polytope complex $P=\bigcup_\ell P_\ell$. \\

In the context of Theorem~\ref{ConeTheorem}, let the monoid $Q_\ell$ denote the image of $\iota_\cC$ in each $\La_\ell$, and set 
$$
Q = \bigcup_\ell Q_\ell = \mathrm{image}(\iota_\cC) \subset \La(\N).
$$
By Theorem~\ref{ConeTheorem}, each $Q_\ell$ is a fine, saturated, sharp, and finite index submonoid of $\La_\ell(\N)=\La_\ell \cap \mathbb{N}^{9g - 9}.$
Therefore, 
$$
D_\ell=\mathrm{Proj}(\mathbb{C}[Q_\ell])
$$
is an irreducible toric variety. The moment polytope of $D_\ell$, however, will be only $\Q$-linearly isomorphic to $P_\ell$;  It may not be exactly $P_\ell$. That's because the generators of free groups $Q_\ell^{gp}$ and $\La_{\ell}^{gp}$ could be different.

\begin{example}
Suppose $M\subset \N^2$ is the set of points $(x,y)$ such that $2\mid(x+y)$. Both $M$ and $\N^2$ are saturated. However, $(\N^2)^{gp}=\Z^2$ in the standard sense, while $M^{gp}$ is generated by $(2,0)$ and $(1,1)$. The projective varieties associated to both $\N^2$ and $M$ are $\P^1$. The inclusion $M\subset \N^2$ implies that the latter $\P^1=\tn{Proj}(\C[M])$ is a $\Z_2$-quotient of $\P^1=\tn{Proj}(\C[\N^2])$. More precisely, 
$$
\tn{Proj}(\C[\N^2])=\P^1[1,1] \quad\tn{and}\quad \tn{Proj}(\C[M])=\P^1[2,2].
$$
\end{example}
\medskip

Theorem~\ref{ThmC} identifies $\mathrm{Sk}_g^{\mathrm{gd}}$ with the Stanley–Reisner ring $\mathbb{C}[Q]$; consequently, $D_g=\partial \ov\cX_g$ decomposes into a union of components:
$$
D_g=\mathrm{Proj}(\mathrm{Sk}_g^{\mathrm{gd}}) = \bigcup_{\ell} D_\ell
$$
in accordance with the polyhedral decomposition of $P$. This correspondence matches the dual intersection complex of the boundary divisor $D_g$  with the dual of the polyhedral complex $P$. Since $P$ is a PL sphere, the dual of the polyhedral complex $P$ is also a PL sphere. This finishes the proof of Theorem~\ref{MainTheorem}. \qed

\section{Proof of Theorem~\ref{f-functions}}\label{sec:DTandf-functions}

We prove Theorem~\ref{f-functions}.(I) by reducing it to a simpler statement, repositioning some of the $U$-curves, and applying a bigon-detection lemma. The proof of Theorem~\ref{f-functions}.(II) is similar but simpler, so we leave it to the reader.\\

Consider the configuration related to Theorem~\ref{f-functions}.(I), illustrated in Figure~\ref{f-functionsFigure}, where $a \in \cP$ belongs to two pairs of pants in the given pants decomposition, and $a_1,a_2,a_3,a_4$ are the other four curves in $\cP$ that belong to the two pairs of pants containing $a$. As in \eqref{shrunken_eq}, we denote the two parallel copies of $a$ that bound the annulus $A$ around $a$ by $a^+$ and $a^-$.  We denote the shrunken pair of pants containing $a_1, a_4,$ and $a_+$ by $F^+$ and the one containing $a_2, a_3,$ and $a^-$ by $F^-$. For every multicurve $\mf{m}$ on $\Si$, the intersection of $\mf{m}$ with $F^+ \cup F^-$ consists of some $U$- and $\beta$-arcs, the number of each of which is determined by $n=i(\mf{m},a)$ and
$$
(n_1,n_2,n_3,n_4) = \lrp{i(\mf{m},a_1),i(\mf{m},a_2),i(\mf{m},a_3),i(\mf{m},a_4)}.
$$
The goal is to express the value of $n' = i(\mf{m},a')$ in terms of the twist coordinate $t$ along $a$ and the tuple $(n,n_1,\ldots,n_4)$.\\

By the {\bf Bigon Criterion} (cf.\ \cite{FLP}), $n'$ is the count of intersection points between two (transversely intersecting) representatives if and only if there is no bigon between them. Moreover, if there is a bigon, there always exists an “innermost” bigon whose interior is disjoint from both $\mf{m}$ and $a'$. Thus, one can remove such a bigon and continue inductively to eliminate all bigons and achieve the minimal intersection number. Having drawn a standard presentation of $\mf{m}$ as in the previous section, it may form some bigons with $a'$ (similarly, with $a''$). However, since $a'$ lies inside $F^+ \cup F^- \cup A$, any innermost bigon between $\mf{m}$ and $a'$ must lie entirely within this region. Otherwise, there would be a bigon between $\mf{m}$ and one of the $a_i$, contradicting the fact that any standard representative intersects the $a_i$ minimally. Therefore, in proving Theorem~\ref{f-functions}.(I), we may restrict our attention to the four-holed sphere $F^+ \cup F^- \cup A$ (similarly, a one-holed torus in part II) and view $\mf{m}$ as a multiarc ending at a prescribed number of points on $a_1,\ldots,a_4$ (similarly, on $a_1$ in part II).\\

With this perspective, recall that the quantities
$$
\aligned
&\de_{14} = \max\left\{ \frac{n_1+n_4 - n}{2}, 0 \right\}, \quad \de_{23} = \max\left\{ \frac{n_2+n_3 - n}{2}, 0 \right\},\\
&u_1 = \max\left\{ \frac{n_1 - n_4 - n}{2}, 0 \right\}, \quad u_4 = \max\left\{ \frac{n_4 - n_1 - n}{2}, 0 \right\},\\
&u_2 = \max\left\{ \frac{n_2 - n_3 - n}{2}, 0 \right\}, \quad u_3 = \max\left\{ \frac{n_3 - n_2 - n}{2}, 0 \right\}
\endaligned
$$
in the statement of Theorem~\ref{f-functions}.(I) count $\beta$- and $U$-arcs in $F^\pm$ and are related to the intersection points with $\mf{m}$ as follows. Regardless of the twist coordinate $t$ along $a$:
\begin{itemize}
\item $\de_{14}$ is the number of $\beta$-arcs between $a_1$ and $a_4$ -- each such arc contributes $1$ to $i(\mf{m},a')$;
\item $\de_{23}$ is the number of $\beta$-arcs between $a_2$ and $a_3$ -- each such arc contributes $1$ to $i(\mf{m},a')$;
\item for $i=1,\ldots,4$, $u_i$ is the number of $U_i$-arcs starting and ending on $a_i$ -- each such $U$-arc contributes $2$ to $i(\mf{m},a')$.
\end{itemize}

\begin{figure}
\begin{pspicture}(-17,-3.5)(0,2.5)
\psset{unit=1.5cm}

\psarc(-6,1.35){.6}{180}{360}
\psarc(-6,-1.85){.6}{0}{180}

\psline(-7.4,1.35)(-7.4,-1.85)
\psline(-4.6,1.35)(-4.6,-1.85)

\psellipticarc[linestyle=dashed, dash=2pt](-6,-.5)(1.4,0.1){0}{180}
\psellipticarc(-6,-.5)(1.4,0.1){180}{360} \rput(-5.3,-.2){\small $a^+$}

\psellipticarc[linestyle=dashed, dash=2pt](-6,0)(1.4,0.1){0}{180}
\psellipticarc(-6,0)(1.4,0.1){180}{360} \rput(-5.3,-.7){\small $a^-$}

\psellipticarc[linecolor=red, linestyle=dashed, dash=2pt](-6,-0.25)(0.2,1){-90}{90} \rput(-6.3,.4){\small $a'$}
\psellipticarc[linecolor=red](-6,-0.25)(0.2,1){90}{270} 

\psellipse(-7,1.35)(0.4,0.1) \rput(-7,1.6){\small $a_1$}
\psellipse(-5,1.35)(0.4,0.1) \rput(-5,1.6){\small $a_4$}
\psellipse(-7,-1.85)(0.4,0.1) \rput(-7,-2.1){\small $a_2$}
\psellipse(-5,-1.85)(0.4,0.1) \rput(-5,-2.1){\small $a_3$}
\end{pspicture}
\caption{Configuration of Theorem~\ref{f-functions}.(I)}
\label{f-functionsFigure}
\end{figure}

Subtracting these $t$-independent contributions from $n'$ and removing the corresponding $U$- and $\beta$-arcs from the picture (which reduces $n_i$ to $\wt{n}_i$), Theorem~\ref{f-functions}.(I) reduces to the following simplified version.

\begin{proposition}\label{f-functions-simple}
(I) With notation as above, and supposing the curve $a$ belongs to two pairs of
pants, if $n \geq n_1 + n_4$ and $n \geq n_2 + n_3$, then
$$
n' = f(n,t,n_1,n_2,n_3,n_4) \quad \tn{and} \quad n'' = f(n,t-n,n_1,n_2,n_3,n_4),
$$
where
\bEq{simpleF}
f(n,t,n_1,n_2,n_3,n_4) = \max \Big\{ \abs{2t + n - n_1 - n_3},\; n_1 + n_3 - n,\; n_2 + n_4 - n \Big\}.
\eEq

(II) Similarly, if $a$ only belongs to one pair of pants in $\cP$ as in Figure~\ref{dual curve}-Right, and $n \geq n_1/2$, then
$$
n' = |t| \quad \tn{and} \quad n'' = |t - n|.
$$
\end{proposition}

\vskip.1in
We begin by outlining the proof of Theorem~\ref{f-functions-simple} to develop some intuition, which we will make precise in the remainder of this section. Throughout the rest of this section, we assume  $n \geq n_1 + n_4$ and $n \geq n_2 + n_3$. Moreover, when drawing a curve in standard position we choose the $n$ points on the front halves $\al^+$ and $\al^-$ of $a^+$ and $a^-$ so that the intersection points of $\mf{m}$ with $a'$ happen in $A$. That means
\begin{itemize}
\item in $F^+$, the $\beta$-arcs connecting $a^+$ and $a_1$ are drawn on the left of $a'$;
\item in $F^+$, the $\beta$-arcs connecting $a^+$ and $a_4$ are drawn on the right of $a'$;
\item in $F^+$, the $U$-arcs connecting $a^+$ to itself are drawn on the right of $a'$;
\item in $F^-$, the $\beta$-arcs connecting $a^-$ and $a_3$ are drawn on the right of $a'$;
\item in $F^-$, the $\beta$-arcs connecting $a^+$ and $a_2$ are drawn on the left of $a'$;
\item in $F^-$, the $U$-arcs connecting $a^-$ to itself are drawn on the left of $a'$;
\end{itemize}
see Figure~\ref{fig:strict-standard}. The ends points are connected in $A$ in the order described by the twist coordinate $t$ along $a$.\\

\begin{figure}
\begin{pspicture}(-17,-3.5)(0,1.5)
\psset{unit=1.4cm}

\psarc(-6,1.35){.6}{180}{360}
\psarc(-6,-1.85){.6}{0}{180}

\psline(-7.4,1.35)(-7.4,-1.85)
\psline(-4.6,1.35)(-4.6,-1.85)

\psellipticarc[linestyle=dashed, dash=2pt](-6,-.5)(1.4,0.1){0}{180}
\psellipticarc(-6,-.5)(1.4,0.1){180}{360} 

\psellipticarc[linestyle=dashed, dash=2pt](-6,0)(1.4,0.1){0}{180}
\psellipticarc(-6,0)(1.4,0.1){180}{360}

\psellipticarc[linecolor=red, linestyle=dashed, dash=2pt](-6,-0.25)(0.2,1){-90}{90} \rput(-6.3,.4){\small $a'$}
\psellipticarc[linecolor=red](-6,-0.25)(0.2,1){90}{270} 

\pscurve[linecolor=blue](-5.5,-.1)(-5.55,0.5)(-5.7,.85)
\pscurve[linecolor=blue, linestyle=dashed](-5.7,.85)(-5.1,.8)(-4.6,.5)
\pscurve[linecolor=blue](-4.6,.5)(-5.1,0.2)(-5.2,-.09)

\pscurve[linecolor=blue](-6.5,-.57)(-6.25,-1.25)(-6.3,-1.35)
\pscurve[linecolor=blue, linestyle=dashed](-6.3,-1.35)(-6.9,-1.4)(-7.4,-1.25)
\pscurve[linecolor=blue](-7.4,-1.25)(-6.8,-.75)(-6.8,-.565)

\psellipse(-7,1.35)(0.4,0.1) \rput(-7,1.6){\small $a_1$} \pscurve[linecolor=blue](-7,1.25)(-6.8,.5)(-6.8,-.05)
\psellipse(-5,1.35)(0.4,0.1) \rput(-5,1.6){\small $a_4$}\pscurve[linecolor=blue](-5,1.25)(-5.2,.5)(-5.3,-.08)
\psellipse(-7,-1.85)(0.4,0.1) \rput(-7,-2.1){\small $a_2$}\pscurve[linecolor=blue](-7,-1.75)(-6.7,-.8)(-6.67,-.58)
\psellipse(-5,-1.85)(0.4,0.1) \rput(-5,-2.1){\small $a_3$}\pscurve[linecolor=blue](-5,-1.75)(-5.3,-.9)(-5.3,-.58)

\end{pspicture}
\caption{Standard presentation of arcs of $\mf{m}$ in $F^+$ and $F^-$.}
\label{fig:strict-standard}
\end{figure}

{\bf Outline of the proof.} There are three cases in \eqref{simpleF}:
\begin{itemize}
\item[(i)] $n \geq n_1 + n_3$ and $n \geq n_2 + n_4$;
\item[(ii)] $n \geq n_1 + n_3$ and $n < n_2 + n_4$.
\item[(iii)] $n < n_1 + n_3$ and $n \geq n_2 + n_4$;
\end{itemize}

Note that, by assumption, we cannot have both $n < n_1 + n_3$ and $n < n_2 + n_4$.\\

Let $U_+$-arcs denote the $U$-arcs starting and ending on $a_+$ in $F^+$, and $U_-$-arcs denote the $U$-arcs starting and ending on $a_-$ in $F^-$.\\

For simplicity, assume $n \!-\!n_1\! -\! n_3$ is even.
In case (i), if $t\leq  \frac{n_1+n_3-n}{2}$, there will be $U_\pm$-curves that form bigons with $a'$. However, we can eliminate these bigons by exchanging $\frac{n-(n_1 + n_3)}{2}$ of the $U_\pm$-arcs with oppositely rotating $U$-arcs, as illustrated in Figure~\ref{Type-i-move}. Figure~\ref{Type-i-move}.Left illustrates the standard presentation of $a'$ itself, which makes two bigons with $a'$. Using the $U$-arc that goes around $a_1$ instead of the one that goes around $a_4$, as shown in Figure~\ref{Type-i-move}.Right,  we obtain a non-standard presentation with no bigon.\\  

This move corresponds to changing the twist coordinate to 
$$
T = t + \frac{n - n_1 - n_3}{2},
$$
so that the configuration corresponding to $t' = 0$ has no intersections with $a'$. Therefore, we obtain
$$
n' = 2|T|.
$$
\begin{figure}
\begin{pspicture}(-13,-3.5)(0,2.5)
\psset{unit=1.4cm}

\psarc(-6,1.35){.6}{180}{360}
\psarc(-6,-1.85){.6}{0}{180}

\psline(-7.4,1.35)(-7.4,-1.85)
\psline(-4.6,1.35)(-4.6,-1.85)

\psellipticarc[linestyle=dashed, dash=2pt](-6,-.5)(1.4,0.1){0}{180}
\psellipticarc(-6,-.5)(1.4,0.1){180}{360} 

\psellipticarc[linestyle=dashed, dash=2pt](-6,0)(1.4,0.1){0}{180}
\psellipticarc(-6,0)(1.4,0.1){180}{360}

\psellipticarc[linecolor=red, linestyle=dashed, dash=2pt](-6,-0.25)(0.2,1){-90}{90} \rput(-6.3,.4){\small $a'$}
\psellipticarc[linecolor=red](-6,-0.25)(0.2,1){90}{270} 

\pscurve[linecolor=blue, linestyle=dashed](-5.7,.85)(-5.1,.8)(-4.6,.5)
\pscurve[linecolor=blue](-4.6,.5)(-5.1,0.2)(-5.2,-.09)(-4.6,-.3)

\pscurve[linecolor=blue, linestyle=dashed](-6.3,-1.35)(-6.9,-1.4)(-7.4,-1.25)
\pscurve[linecolor=blue](-7.4,-1.25)(-6.8,-.75)(-6.8,-.565)(-7.4,-.3)


\pscurve[linecolor=blue](-6.3,-1.35)(-6.25,-1.25)(-6.5,-.57)(-5.5,-.1)(-5.55,0.5)(-5.7,.85)

\psellipticarc[linecolor=blue, linestyle=dashed, dash=4pt](-6,-.3)(1.4,0.1){0}{180}

\psellipse(-7,1.35)(0.4,0.1) \rput(-7,1.6){\small $a_1$}
\psellipse(-5,1.35)(0.4,0.1) \rput(-5,1.6){\small $a_4$}
\psellipse(-7,-1.85)(0.4,0.1) \rput(-7,-2.1){\small $a_2$}
\psellipse(-5,-1.85)(0.4,0.1) \rput(-5,-2.1){\small $a_3$}

\psarc(-1,1.35){.6}{180}{360}
\psarc(-1,-1.85){.6}{0}{180}

\psline(-2.4,1.35)(-2.4,-1.85)
\psline(.4,1.35)(.4,-1.85)

\psellipticarc[linestyle=dashed, dash=2pt](-1,-.5)(1.4,0.1){0}{180}
\psellipticarc(-1,-.5)(1.4,0.1){180}{360} 

\psellipticarc[linestyle=dashed, dash=2pt](-1,0)(1.4,0.1){0}{180}
\psellipticarc(-1,0)(1.4,0.1){180}{360}

\psellipticarc[linecolor=red, linestyle=dashed, dash=2pt](-1,-0.25)(0.2,1){-90}{90} \rput(-1.3,.4){\small $a'$}
\psellipticarc[linecolor=red](-1,-0.25)(0.2,1){90}{270} 

\psellipse(-2,1.35)(0.4,0.1) \rput(-2,1.6){\small $a_1$}
\psellipse(0,1.35)(0.4,0.1) \rput(0,1.6){\small $a_4$}
\psellipse(-2,-1.85)(0.4,0.1) \rput(-2,-2.1){\small $a_2$}
\psellipse(0,-1.85)(0.4,0.1) \rput(0,-2.1){\small $a_3$}

\pscurve[linecolor=blue](-1.3,-1.35)(-1.6,-.6)(-1.6,-.1)(-1.55,0.5)(-1.3,.85)
\pscurve[linecolor=blue, linestyle=dashed](-1.3,.85)(-1.45,0.9)(-2.4,.5)
\pscurve[linecolor=blue](-2.4,.5)(-2,.3)(-1.9,-.085)(-1.9,-.585)(-2.4,-1.25)

\pscurve[linecolor=blue, linestyle=dashed](-1.3,-1.35)(-1.9,-1.4)(-2.4,-1.25)



\end{pspicture}
\caption{The bigon removing move in case (i) for $n_i=0$, $n=2$, and $t=-1$ ($\mf{m}=a'$). Left: $\mf{m}$ drawn in standard position that makes two bigons with $a'$. Right: after moving the $U_+$-arc to the opposite side, there is no bigon.}
\label{Type-i-move}
\end{figure}

In case (ii), move (i) can be done to partially remove some of the bigons and the intersection number with $a'$ remains fixed within some range; it will grow linearly outside of that.\\  

In case (iii), there are no bigon but changing $t$ within the range $0\leq t\leq \frac{ n_1 +n_3-n}{2}$ trades one intersection with $a'$ in front for an intersection with $a'$ in back so that the total intersection remains the fixed number $n_1+n_3-n$. Outside this range, $n'$ will grow linearly with respect to $\pm t$. \\  

To make this precise, we draw $\mf{m}$ in a semi-standard position and prove a criterion for when a semi-standard presentation makes a bigon with $a'$.

\begin{definition}\label{ss-pres}
Given $(n,n_1,\ldots,n_4)\in \N^5$ satisfying $n \geq n_1 + n_4$ and $n \geq n_2 + n_3$ and $T\in \Z$, a {\bf semi-standard} presentation of a multiarc $\mf{m}$ in $F^-\cup A \cup F^+$ consists of 
\begin{itemize}
\item $n_1$ $\beta$-arcs in $F^+$ on the left of $a'$ connecting $a^+$ and $a_1$;
\item $n_4$ $\beta$-arcs in $F^+$ on the right of $a'$ connecting $a^+$ and $a_4$;
\item $n_2$ $\beta$-arcs in $F^-$ on the left of $a'$ connecting $a^-$ and $a_2$;
\item $n_3$ $\beta$-arcs in $F^-$ on the right of $a'$ connecting $a^-$ and $a_3$;
\item $x_{1}$ $U$-type arcs in $F^+$  on the left of $a'$ rotating around $a_1$ instead, connecting $a^+$ to itself;
\item $x_{2}$ standard $U$-arcs in $F^-$ on the left of $a'$ connecting $a^-$ to itself;
\item $x_{3}$ $U$-type arcs in $F^-$ on the right of $a'$ rotating around $a_3$ instead, connecting $a^-$ to itself;
\item $x_{4}$ standard $U$-arcs in $F^+$ on the right of $a'$ connecting $a^+$ to itself ;
\end{itemize} 
such that 
$$
(n_4+2x_{4})+(n_1+2x_{1})=n= (n_2+2x_{2})+(n_3+2x_{3});
$$
see Figure~\ref{Semi-std_fig}.
Furthermore, the $n$ end-points on the selected halves $\al^\pm$ (of $a^\pm$) are connected in $A$ using the twist coordinate $T$ in the standard way. If $n=0$, then $\mf{m}$ is just $T$ parallel copies of $a$ in $A$.
\end{definition}

\begin{figure}
\begin{pspicture}(-17,-3.5)(0,1.5)
\psset{unit=1.4cm}

\psarc(-6,1.35){.6}{180}{360}
\psarc(-6,-1.85){.6}{0}{180}

\psline(-7.4,1.35)(-7.4,-1.85)
\psline(-4.6,1.35)(-4.6,-1.85)

\psellipticarc[linestyle=dashed, dash=2pt](-6,-.5)(1.4,0.1){0}{180}
\psellipticarc(-6,-.5)(1.4,0.1){180}{360} 

\psellipticarc[linestyle=dashed, dash=2pt](-6,0)(1.4,0.1){0}{180}
\psellipticarc(-6,0)(1.4,0.1){180}{360}

\psellipticarc[linecolor=red, linestyle=dashed, dash=2pt](-6,-0.25)(0.2,1){-90}{90} \rput(-6.3,.4){\small $a'$}
\psellipticarc[linecolor=red](-6,-0.25)(0.2,1){90}{270} 

\pscurve[linecolor=blue](-5.5,-.1)(-5.55,0.5)(-5.7,.85)
\pscurve[linecolor=blue, linestyle=dashed](-5.7,.85)(-5.1,.8)(-4.6,.5)
\pscurve[linecolor=blue](-4.6,.5)(-5.1,0.2)(-5.2,-.09)
 \rput(-4.8,.8){\small $x_4$}
 
\pscurve[linecolor=blue](-6.5,-.1)(-6.45,0.5)(-6.3,.85)
\pscurve[linecolor=blue, linestyle=dashed](-6.3,.85)(-6.9,.8)(-7.4,.5)
\pscurve[linecolor=blue](-7.4,.5)(-7,0.2)(-7,-.07)
 \rput(-7.2,.8){\small $x_1$}
\pscurve[linecolor=blue](-6.5,-.57)(-6.25,-1.25)(-6.3,-1.35)
\pscurve[linecolor=blue, linestyle=dashed](-6.3,-1.35)(-6.9,-1.4)(-7.4,-1.25)
\pscurve[linecolor=blue](-7.4,-1.25)(-6.8,-.75)(-6.8,-.565)
 \rput(-7.2,-.9){\small $x_2$}
 
\pscurve[linecolor=blue](-5.5,-.57)(-5.75,-1.25)(-5.7,-1.35)
\pscurve[linecolor=blue, linestyle=dashed](-5.7,-1.35)(-5.1,-1.4)(-4.6,-1.25)
\pscurve[linecolor=blue](-4.6,-1.25)(-5.2,-.75)(-5.2,-.565)
 \rput(-4.8,-.9){\small $x_3$}

\psellipse(-7,1.35)(0.4,0.1) \rput(-7,1.6){\small $a_1$} \pscurve[linecolor=blue](-7,1.25)(-6.8,.5)(-6.8,-.05)
\psellipse(-5,1.35)(0.4,0.1) \rput(-5,1.6){\small $a_4$}\pscurve[linecolor=blue](-5,1.25)(-5.2,.5)(-5.3,-.08)
\psellipse(-7,-1.85)(0.4,0.1) \rput(-7,-2.1){\small $a_2$}\pscurve[linecolor=blue](-7,-1.75)(-6.7,-.8)(-6.67,-.58)
\psellipse(-5,-1.85)(0.4,0.1) \rput(-5,-2.1){\small $a_3$}\pscurve[linecolor=blue](-5,-1.75)(-5.3,-.9)(-5.3,-.58)

\end{pspicture}
\caption{Semi-standard presentation of arcs of $\mf{m}$ in $F^+$ and $F^-$.}
\label{Semi-std_fig}
\end{figure}

\begin{remark}
The twist coordinate $t$ of the standard presentation of $\mf{m}$ and the twist coordinate $T$ of a semi-standard presentation of $\mf{m}$ are related by 
$$
t=T-x_{1}-x_{3},
$$
because in order to modify a  semi-standard presentation to the standard one we need to rotate $x_{1}$ points on $a^+$ and $x_{3}$ points on $a^-$ backward to the other side of $a'$. Also, this process does not create any bigons with $a$. Therefore, just like the standard presentation, any semi-standard presentation of $\mf{m}$ intersects $a$ minimally. For $T=0$, the number of intersection points between $\mf{m}$ and $a'$ is $|(n_2+2x_2)-(n_1+2x_1)|=|(n-n_1-n_3)-2(x_1+x_3)|$.
\end{remark}

\begin{lemma}\label{Bigon-detection}
With notation as in Definition~\ref{ss-pres}, $\mf{m}$ makes a bigon with $a'$ if and only if 
\begin{enumerate}
\item $x_4>0$ and $ \de< T < 0$;
\item $x_1>0$ and $  \de> T > 0$;
\item $x_3>0$ and $  \de> T > 0$;
\item $x_2>0$ and $  \de< T < 0$;
\end{enumerate}
where 
$$
\de=(n_1+2x_1)-(n_2+2x_2)=(n_3+2x_3)-(n_4+2x_4).
$$
\end{lemma}

\begin{proof}
By Definition~\ref{ss-pres}, all the intersection points between $\mf{m}$ and $a'$ are in $A$. The intersection of $a'$ and $A$ consists of two (vertical) arcs. We call the one in $N(\Gamma)$ the {\bf front arc}, as $N(\Gamma)$ is always the front half of the pictures we draw. Suppose $\mf{m}$ and $a'$ make a bigon.  The two vertices of any such bigon can not both be on the same arc, say $p_+$ will be on the front arc and the other, say $p_-$ will be on the backside arc. That is because any bigon with two end vertices on the same arc would include a sub-bigon between $\mf{m}$ and $a$, but we know that $\mf{m}$ and $a$ intersect minimally.\\

There are four symmetric possibilities for the innermost bigon, depending on the location of the side of that, say $m$, which is part of $\mf{m}$:
\begin{enumerate}
\item $m$ starts from $p_+$, intersects $\al^+$ to the right of $a'$ at the first point after $a'$, followed by a $U$-curve turning around $a_4$ intersecting $\al^+$ again at the last point, and finally intersecting $a'$ again at $p_-$ on the back via a line segment in $A$ that connects this last point on $\al^+$ to some point on $\al^-$ via some negative twist $T$. The reason the intersection points with $\al^+$ are the first and the last points after $a'$ is because we are considering the innermost bigon.
\item $m$ starts from $p_+$, intersects $\al^+$ to the left of $a'$ at the lest point before $a'$, followed by a $U$-curve turning around $a_1$ intersecting $\al^+$ again at the first point, and finally intersecting $a'$ again at $p_-$ on the back via a line segment in $A$ that connects this first point on $\al^+$ to some point on $\al^-$ via some positive twist $T$.
\item $m$ starts from $p_+$, intersects $\al^-$ to the right of $a'$ at the first point after $a'$, followed by a $U$-curve turning around $a_3$ intersecting $\al^-$ again at the last point, and finally intersecting $a'$ again at $p_-$ on the back via a line segment in $A$ that connects this last point on $\al^-$ to some point on $\al^+$ via some positive twist $T$.
\item $m$ starts from $p_+$, intersects $\al^-$ to the right of $a'$ at the last point before $a'$, followed by a $U$-curve turning around $a_2$ intersecting $\al^-$ again at the first point, and finally intersecting $a'$ again at $p_-$ on the back via a line segment in $A$ that connects this last point on $\al^-$ to some point on $\al^+$ via some negative twist $T$.
\end{enumerate}
The cases 1-4  correspond to cases 1-4 of the statement, respectively. 
\end{proof}

{\bf Proof of Proposition~\ref{f-functions-simple}.} We include the details of the three cases (i)-(iii) outlined above.
For a curve drawn in the standard positive there are $n_1$ intersection points with $\al^+$  and $n-n_3$ intersection points with $\al^-$ on the left of $a^-$. So their difference is 
$$
n-(n_1+n_3)\geq 0.
$$
In case (i), comparing  twice the number of $U$ arcs in the standard position, that is
$$
(n-(n_1+n_4))+(n-(n_2+n_3)),
$$
with $n-(n_1+n_3)$, since $n-(n_2+n_4)\geq 0$, we have 
$$
(n-(n_1+n_4))+(n-(n_2+n_3))\geq n-(n_1+n_3). 
$$
Therefore, if $n-(n_1+n_3)$ is even, there is a semi-standard presentation with $\de=0$ and 
$$
x_1+x_3= (n-(n_1+n_3))/2.
$$
Then, by Lemma~\ref{Bigon-detection},  for any $T$, there is no bigon between such semi-standard presentation and $a'$. A simple count of intersection points in $A$ yields 
$$
n'=2|T|= 2|t+x_1+x_3|=|2t+n-n_1-n_3|.
$$

If $n-(n_1+n_3)$ is odd, there is a semi-standard presentation with $\de=\pm 1$. Again, by Lemma~\ref{Bigon-detection},  for any $T$, there is no bigon between such semi-standard presentation and $a'$. A simple count of intersection points in $A$ yields 
$$
n'=|2t+n-n_1-n_3|.
$$

In case (ii), moving all the $U$-arcs in $F^+$ to the $a_1$-side and all the $U$-arcs in $F^-$ to the $a_3$-side we get a semi-standard presentation with 
$$
\de=n-n_4-n_2<0,~x_4,x_2=0.
$$
Again, by Lemma~\ref{Bigon-detection},  for any $T$, there is no bigon between such semi-standard presentation and $a'$. A simple count of intersection points in $A$ yields 
$$
n'=\tn{max}\lrc{|2t+n-n_1-n_3|, n_2+n_4-n}.
$$

Finally, in case (iii), the standard presentation satisfies 
$$
x_1,x_3=0\quad \tn{and}\quad \de=n_1+n_3-n>0.
$$ 
By Lemma~\ref{Bigon-detection},  for any $T$, there is no bigon between such semi-standard presentation and $a'$. A simple count of intersection points in $A$ yields 
$$
n'=\tn{max}\lrc{|2t+n-n_1-n_3|, n_1+n_3-n}.
$$
This establishes the formula for $n'$ in \eqref{simpleF}. The formula for $n''$ immediately follows from the fact that $a''$ is a Dehn Twist of $a'$ along $a$.\\

This finishes the proof of Proposition~\ref{f-functions-simple}.(I) and thus Theorem~\ref{f-functions}.(I).
Proof of Proposition~\ref{f-functions-simple}.(II) is similar (and simpler).  Thus we won't repeat that and leave it to the reader.
\qed\\

\section{Proof of Theorem~\ref{ThmC}}\label{Sec:proof-of-MainTheorem}

In order to prove Theorem~\ref{ThmC}, we will use a modified version of DT coordinates due to \cite{BKL}.\\

For every multicurve $\mf{m}$, let 
$$
\nu(\mf{m})=(n_1,\ldots,n_{3g-3},t_1,\ldots,t_{3g-3}) \in \N^{3g-3}\times \Z^{3g-3}
$$
denote its DT coordinates, and let
$$
\iota_{\cC}(\mf{m})=(n_1,\ldots,n_{3g-3},n'_1,\ldots,n'_{3g-3}, n''_1,\ldots,n''_{3g-3})\in \N^{9g-9}
$$
denote its intersection coordinates with the curves in the collection $\cC=(\cP,\cQ,\cQ_+)$. \\

For every pair of multicurves $\mf{m}$ and $\mf{m}'$, recall that 
\bEq{decomposition2}
\mf{m}\ast\mf{m}'=\sum_{s} \la_s \mf{m}_s,
\eEq
where $s$ runs over all possible smoothings of $\mf{m}\cup \mf{m}'$ (in general position), and the coefficient $\la_s$ collects the contributions of trivial components (a trivial curve contributes $-2$) appearing in a smoothing. \\ 

It is clear that for every curve $c$ and every smoothing $\mf{m}_s$ in \eqref{decomposition2}, we have 
$$
i(\mf{m}_s,c)\leq i(\mf{m},c)+ i(\mf{m}',c),
$$
with equality if and only if the smoothing $\mf{m}_s$ does not form any bigon with $c$ (assuming $\mf{m}$ and $\mf{m}'$ were originally positioned to avoid any bigons with $c$). \\

With the degree of each multicurve defined as 
$$
 |\mf{m}|=\sum_{c\in \cC} i(\mf{m},c),
$$
we conclude that  
\begin{equation}\label{norm-inequality}
|\mf{m}_s|\leq |\mf{m}|+|\mf{m}'|\qquad \forall~s,
\end{equation}
with equality if and only if 
$$
\iota_{\cC}(\mf{m}_s)= \iota_{\cC}(\mf{m})+\iota_{\cC}(\mf{m}')\in \N^{9g-9}.
$$
Theorem~\ref{ThmC} states that this equality holds if and only if $\mf{m}$ and $\mf{m}'$ lie in the same cone $\La_\ell$ for some $\ell$, and $\mf{m}_s$ is a uniquely determined geometric sum $m_\#$ of the two multicurves $\mf{m}$ and $\mf{m}'$. \\

\begin{figure}
    \begin{pspicture}(-10,-1.5)(0,1.5)
    \psset{unit=1cm}

    \psellipse(-7,1)(0.4,0.1) \rput(-7,1.3){\small $a_3$}
    \psellipse(-5,1)(0.4,0.1) \rput(-5,1.3){\small $a_2$}
    \psellipticarc[linestyle=dashed, dash=1pt](-6,-1)(.4,.1){0}{180}\rput(-6,-1.3){\small $a_1$}
    \psellipticarc(-6,-1)(.4,.1){180}{360}\rput(-6,-1.3)

    \psarc(-6,1){.6}{180}{360}
    \psbezier(-4.6,1)(-4.6,-.5)(-5.6,.5)(-5.6,-1)
    \psbezier(-7.4,1)(-7.4,-.5)(-6.4,.5)(-6.4,-1)

    \pscurve[linecolor=red](-6,-1.1)(-6.1,0)(-6,.4)
    \pscurve[linecolor=red,linestyle=dashed, dash=1pt](-6,.4)(-5.6,0.1)(-5.5,-.4)\rput(-5.9,-.3){\small $\mf{m}$}
    \pscurve[linecolor=red](-5.5,-.4)(-5.7,-.7)(-5.8,-1.05)

    \psarc[linecolor=blue](-6,1){.9}{185}{355}\rput(-6.7,0.3){\small $\mf{m}'$}

		\psline{->}(-4.5,-.3)(-3,-.3)
		
     \psellipse(-2,1)(0.4,0.1) \rput(-2,1.3){\small $a_3$}
    \psellipse(0,1)(0.4,0.1) \rput(0,1.3){\small $a_2$}
    \psellipticarc[linestyle=dashed, dash=1pt](-1,-1)(.4,.1){0}{180}\rput(-1,-1.3){\small $a_1$}
    \psellipticarc(-1,-1)(.4,.1){180}{360}

    \psarc(-1,1){.6}{180}{360}
    \psbezier(.4,1)(.4,-.5)(-.6,.5)(-.6,-1) 
    \psbezier(-2.4,1)(-2.4,-.5)(-1.4,.5)(-1.4,-1)

    \pscurve[linestyle=dashed, dash=1pt](-1,.4)(-.6,0.1)(-0.5,-.4)
    \pscurve(-0.5,-.4)(-0.7,-.7)(-0.8,-1.05)
    \psarc(-1,1){.9}{285}{355}
     \pscurve(-1,.4)(-.9,.15)(-.7,0.15)

    \psarc(-1,1){.9}{185}{250}
    \pscurve(-1,-1.1)(-1.1,0)(-1.32,.16)
\rput(1,-.3){\small $\mf{m}_{s_1}~~ +~~ \mf{m}_{s_2}$}


     \psellipse(2,1)(0.4,0.1) \rput(2,1.3){\small $a_3$}
    \psellipse(4,1)(0.4,0.1) \rput(4,1.3){\small $a_2$}
    \psellipticarc[linestyle=dashed, dash=1pt](3,-1)(.4,.1){0}{180}\rput(3,-1.3){\small $a_1$}
    \psellipticarc(3,-1)(.4,.1){180}{360}

    \psarc(3,1){.6}{180}{360}
    \psbezier(4.4,1)(4.4,-.5)(3.4,.5)(3.4,-1) 
    \psbezier(1.6,1)(1.6,-.5)(2.6,.5)(2.6,-1)

    \pscurve[linestyle=dashed, dash=1pt](3,.4)(3.4,0.1)(3.5,-.4)
    \pscurve(3.5,-.4)(3.3,-.7)(3.2,-1.05)
    \psarc(3,1){.9}{285}{355}
        \pscurve(3,-1.1)(3,0)(3.3,0.15)

     \pscurve(3,.4)(2.9,.15)(2.68,.16)
        \psarc(3,1){.9}{185}{250}

    \end{pspicture}
    \caption{$\mf{m}$, $\mf{m}'$ and the two smoothings in their product}
    \label{curve-mult}
    \end{figure}
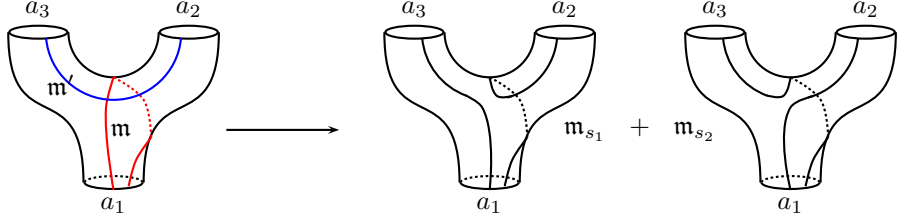

The map $\mf{m} \mapsto \nu(\mf{m})$, given by DT coordinates, does not behave well under the product in the sense above. The vectors $\nu(\mf{m}_s)$ on the right-hand side of \eqref{decomposition} can behave somewhat erratically if there are $U$-arcs in the standard presentations of $\mf{m}$ and $\mf{m}'$. For instance, in the local Figure~\ref{curve-mult}, the $(n,t)$-coordinates of $\mf{m}$ and $\mf{m'}$ are 
$$
\nu(\mf{m})=\big((2,0,0),(0,0,0)\big)\quad \tn{and}\quad \nu(\mf{m}')=\big((0,1,1),(0,0,0)\big)
$$ 
while the $(n,t)$-coordinates of neither of their two smoothings is equal to the sum of these vectors. More precisely, we have 
$$
\nu(\mf{m}_{s_1})=\big((2,1,1),(0,1,0)\big)\quad \tn{and}\quad \nu(\mf{m}_{s_2})=\big((2,1,1),(1,0,-1)\big)
$$
and both have non-trivial twist coordinates.\\

\begin{figure}
\begin{pspicture}(-14,-2.5)(0,2)
\psset{unit=1cm}

\psarc(-6,1.35){.6}{180}{360}
\psarc(-6,-1.85){.6}{0}{180}

\psline(-7.4,1.35)(-7.4,-1.85)
\psline(-4.6,1.35)(-4.6,-1.85)

\psellipticarc[linestyle=dashed, dash=2pt](-6,-0.25)(1.4,0.2){0}{180}
\psellipticarc(-6,-.25)(1.4,0.2){180}{360} \rput(-6,-.7){\small $a$}

\psellipse(-7,1.35)(0.4,0.1) \rput(-7,1.6){\small $a_1$}
\psellipse(-5,1.35)(0.4,0.1) \rput(-5,1.6){\small $a_4$}
\psellipse(-7,-1.85)(0.4,0.1) \rput(-7,-2.1){\small $a_2$}
\psellipse(-5,-1.85)(0.4,0.1) \rput(-5,-2.1){\small $a_3$}
\end{pspicture}
\caption{Two pants including $a$}
\label{two-s-pants}
\end{figure}
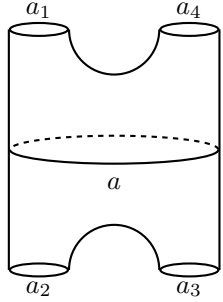

To resolve this issue, in \cite{BKL}, the authors consider a modified version of DT coordinates. They adjust the twist coordinate $t_i$ along $a_i$ as follows:
$$
t'_i= t_i + \# \tn{ number of } U-\tn{arcs turning around } a_i \tn{ on each side}.
$$
In other words, the $t$ coordinate of the pants curve $a$ in Figure~\ref{two-s-pants} is changed to
\begin{equation}\label{MDT}
t'= t+ \tn{max}\lrc{(n_1-(n+n_4))/2,0}+\tn{max}\lrc{ (n_3-(n+n_2))/2,0}.
\end{equation}

We will denote the {\bf modified DT} coordinate vector by $\nu'(\mf{m})$. Then, by \cite[Theorem~9.1]{FKL23}, the product $\mf{m} \ast \mf{m'}$ has a unique {\bf leading order term} $\mf{m}_\#$ satisfying 
\begin{equation}\label{leading-term}
\nu'(\mf{m}_\#)= \nu'(\mf{m})+\nu'(\mf{m}')\in \N^{3g-3}\times \Z^{3g-3}
\end{equation}
and the remaining terms $\mf{m}_s$ in \eqref{decomposition} are strictly smaller than $\mf{m}_\#$ with respect to a certain partial order. We will not need the details of this ordering. For instance, in the example above, ignoring contributions from adjacent pairs of pants, we have
$$
\aligned
&\nu'(\mf{m})=\big((2,0,0),(0,1,0)\big),\quad \nu'(\mf{m}')=\big((0,1,1),(0,0,0)\big),\\
&\nu'(\mf{m}_{s_1})=\big((2,1,1),(0,1,0)\big),\quad \nu'(\mf{m}_{s_2})=\big((2,1,1),(1,0,-1)\big).
\endaligned
$$
Therefore, $\mf{m}_\#=\mf{m}_{s_1}$.

\subsection*{Proof of Theorem~\ref{ThmC}}
\begin{lemma}[Proof of Theorem~\ref{ThmC} - Part 1]\label{Part1}
For each $\La_\ell$, and every two multicurves with $\iota_\cC(\mf{m}),\iota_\cC(\mf{m}')\in\La_\ell$,  the product $\mf{m} \ast \mf{m'}$ has a unique leading order term $\mf{m}_\#$ with  $\iota_\cC(\mf{m}_\#)\in \La_\ell$ satisfying 
$$
|\mf{m}_\#|= |\mf{m}|+|\mf{m}'|.
$$
\end{lemma}
\begin{proof}

Note that the modified DT coordinates are piecewise linear functions of the standard DT coordinates. More precisely, by \eqref{MDT}, the regions in which the transformation is linear are determined by the inequalities  
$$
n_i\geq n_j+n_k~~~\tn{or}~~~n_i\leq n_j+n_k
$$
in each pair of pants $F(a_i,a_j,a_k)$. These inequalities are a subset of the inequalities that determine the regions of linearity $\La'_\ell$ of $\iota_\cC$ in $\overline{\mc{D}_g}$. We conclude that:\\

{\it Over each $\La'_\ell$, the map from DT coordinates to modified DT coordinates is the restriction to $\La'_\ell$ of a rational linear isomorphism $\Psi_\ell\colon \R^{6g-6}\lra \R^{6g-6}$.}\\

Before we use the last observation to prove Lemma~\ref{Part1}, we need to modify the definition of leading term in \eqref{leading-term} in some special cases. Fixing one of the cones $\La'_\ell\subset \ov{\mc{D}_g}$, in Convention~\ref{new-convention}, we modified the DT coordinates of multicurves $\mf{m}$ with $n_i=0$ to have twist coordinates $t_i\leq 0$ instead of $t_i\geq 0$, whenever $n''_i>n'_i$. That's in order for $\iota_\cC$ to define a continuous map between two closed cones $\La'_\ell$ and $\La_\ell$. A similar change must be applied to the definition of leading term in \eqref{leading-term} when such $\mf{m}$ is multiplied with another multicurve $\mf{m}'$ in $\La'_\ell$ as follows.\\

If $n_i(\mf{m})=0$, $\mf{m}= a_i^{\tau}~\mf{m}''$ for some $\tau\geq 0$ such that both $n_i$ and $t_i$ of $\mf{m}''$ are zero. Since the argument below only concerns the $i$-th twist coordinate, we may assume $\mf{m}= a_i^{\tau}$. Then, \eqref{leading-term} states that there is a unique leading term $\mf{m}_\#$ in $a_i^{\tau}\ast \mf{m'}$ such that 
$$
\nu'(\mf{m}_\#)=\nu'(\mf{m}')+ ((0,\ldots,0), (0,\ldots,0, \tau, 0, \ldots,0);
$$
i.e. all the modified DT (as well as standard DT) coordinates of $\mf{m}_\#$ are the same as in $\mf{m}'$; only the $i$-th twist coordinate increases by $\tau$. Topologically, this leading order term is obtained as follows. The standard representative of the multicurve $\mf{m}'$ and $\tau$ copies of $a_i$ intersect at $n_i\cdot \tau$ points in the annulus $A_i$ around $a_i$; see Figure~\ref{Grid}-Left. To get $\mf{m}_\#$, all these intersection points must be resolved positively as in Figure~\ref{Grid}-Right. 
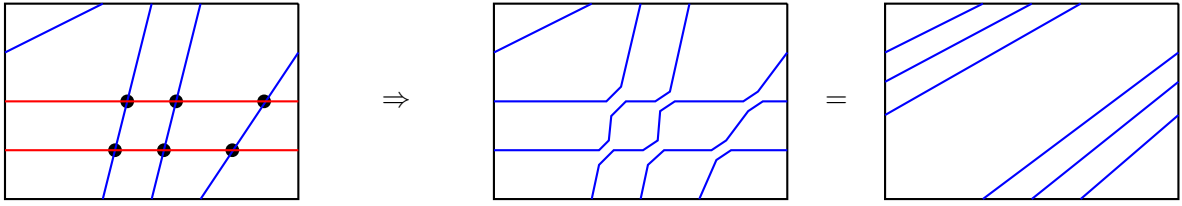
\begin{figure}
\begin{pspicture}(-1,-1.5)(0,1.5)
\psset{unit=1.3cm}

\psline(2.5,1)(2.5,-1)(-.5,-1)(-.5,1)(2.5,1)
\pscircle*(.75,0){.07}\pscircle*(1.25,0){.07}\pscircle*(2.15,0){.07}
\pscircle*(.625,-.5){.07}\pscircle*(1.125,-.5){.07}\pscircle*(1.825,-.5){.07}
\psline[linecolor=blue](.5,-1)(1,1)
\psline[linecolor=blue](1,-1)(1.5,1)
\psline[linecolor=blue](1.5,-1)(2.5,.5)
\psline[linecolor=blue](-.5,.5)(.5,1)

\psline[linecolor=red](-.5,-.5)(2.5,-.5)
\psline[linecolor=red](-.5,0)(2.5,0)

    \rput(3.5,0){$\Rightarrow$}
\psline(7.5,1)(7.5,-1)(4.5,-1)(4.5,1)(7.5,1)

 \psline[linecolor=blue](4.5,0)(5.65,0)(5.8,.15)(6,1)
 \psline[linecolor=blue](4.5,-.5)(5.575,-.5)(5.675,-.4)(5.7,-.15)(5.85,0)(6.15,0)(6.3,0.1)(6.5,1)
 \psline[linecolor=blue](5.5,-1)(5.575,-.65)(5.725,-.5)(6.025,-.5)(6.175,-.4)(6.2,-.1)(6.35,0)(7.05,0)(7.2,0.1)(7.5,.5)
\psline[linecolor=blue](4.5,.5)(5.5,1)
  \psline[linecolor=blue](6,-1)(6.075,-.65)(6.24,-.5)(6.725,-.5)(6.875,-.4)(7.1,-.1)(7.25,0)(7.5,0)
    \psline[linecolor=blue](6.6,-1)(6.775,-.6)(6.925,-.5)(7.5,-.5)
    
    \rput(8,0){$=$}
    
    \psline(11.5,1)(11.5,-1)(8.5,-1)(8.5,1)(11.5,1)
     \psline[linecolor=blue](9.5,-1)(11.5,.5)
\psline[linecolor=blue](10,-1)(11.5,0.2)
\psline[linecolor=blue](10.5,-1)(11.5,-.14)

\psline[linecolor=blue](8.5,.5)(9.5,1)
\psline[linecolor=blue](8.5,0.2)(10,1)
\psline[linecolor=blue](8.5,-.14)(10.5,1)

\end{pspicture}
\caption{Left: The intersection points of $a_i^2$ (red curves) and some $\mf{m'}$ with $n_i=3$ and $t_i=1$ (blue curves) in the annulus $A_i$, presented as a rectangle with the vertical sides identified. Right: the positive resolution that increases $t_i$ by $\tau$.}
\label{Grid}
\end{figure}
However, if we define the $i$-th twist coordinate of $a_i^\tau$ to be $-\tau$, then in order to preserve the leading term property of \eqref{leading-term}, we must instead resolve all the intersection points negatively, as in Figure~\ref{negative-res}. 
\begin{figure}
\begin{pspicture}(3,-1.5)(3,1.5)
\psset{unit=1.3cm}

\psline(7.5,1)(7.5,-1)(4.5,-1)(4.5,1)(7.5,1)

\psline[linecolor=blue](4.5,0)(5.65,0)(5.7,-.15)(5.675,-.4)(5.725,-.5)(6.025,-.5)(6.075,-.65)(6,-1)
\psline[linecolor=blue](6,1)(5.8,.15)(5.85,0)(6.15,0)(6.2,-.1)(6.175,-.4)(6.24,-.5)(6.725,-.5)(6.775,-.6)(6.6,-1)
\psline[linecolor=blue](6.5,1)(6.3,0.1)(6.35,0)(7.05,0)(7.1,-.1)(6.875,-.4)(6.925,-.5)(7.5,-.5)
\psline[linecolor=blue](4.5,-.5)(5.575,-.5)(5.575,-.65)(5.5,-1)
\psline[linecolor=blue](4.5,.5)(5.5,1)
\psline[linecolor=blue](7.5,.5)(7.2,0.1)(7.25,0)(7.5,0)

 \rput(8,0){$=$}
 
     \psline(11.5,1)(11.5,-1)(8.5,-1)(8.5,1)(11.5,1)
\psline[linecolor=blue](10,-1)(9.5,1)
\psline[linecolor=blue](10.5,-1)(10,1)
\psline[linecolor=blue](9.5,-1)(8.5,-.5)
 \psline[linecolor=blue](10.5,1)(11.5,-.5)
 
\end{pspicture}
\caption{Negative resolution of Figure~\ref{Grid}.left, decreasing the twist coordinate from $t_i=1$ to $t_i-\tau=-1$.}
\label{negative-res}
\end{figure}
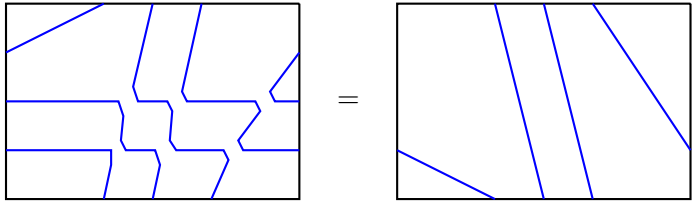
Then, this differently chosen lead term satisfies 
$$
\nu'(\mf{m}_\#)=\nu'(\mf{m}')+ ((0,\ldots,0), (0,\ldots,0, -\tau, 0, \ldots,0);
$$
We conclude the following:\vskip.05in

{\it For every $\La'_\ell$ and every two multicurves $\mf{m},\mf{m'}\in \La'_\ell$  the product $\mf{m} \ast \mf{m'}$ has a unique $\ell$-dependent {\bf leading order term} $\mf{m}_{\ell,\#}$ satisfying 
\begin{equation}\label{leading-term2}
\nu'(\mf{m}_{\ell,\#})= \nu'(\mf{m})+\nu'(\mf{m}')\in  \Z^{6g-6}
\end{equation}
where the twist coordinates for when $n_i=0$ are determined by Convention~\ref{new-convention}.}\\

Now, with this modification, the map $\iota_{\cC}\circ \Psi_\ell^{-1}$ from the coordinate domain $\La_\ell'=\Psi_\ell(\La'_\ell)$ of twisted DT coordinates on multicurves to $\La_\ell$ is linear. Therefore, it follows from (\ref{leading-term2}) that  
$$
\iota_\cC(\mf{m}_{\ell,\#})= \iota_\cC(\mf{m})+\iota_\cC(\mf{m}'),
$$  
for every pair of multicurves with $\iota_\cC(\mf{m}),\iota_\cC(\mf{m}')\in\La_\ell$. Furthermore, $\mf{m}_{\ell,\#}$ is unique. This completes the proof of the lemma.
\end{proof}

{\bf Proof of Theorem~\ref{ThmC} - Part 2.} Suppose $\iota_\cC(\mf{m})$ and $\iota_\cC(\mf{m'})$ do not belong to the same cone $\La_\ell$, but 
$$
\mf{m}~\ov\ast~\mf{m}'\neq 0.
$$
This means there is a multicurve smoothing $m_\#$ of their union satisfying 
\begin{equation}\label{mforlemma}
\iota_\cC(\mf{m}_{\#})= \iota_\cC(\mf{m})+\iota_\cC(\mf{m}');
\end{equation}
see \eqref{norm-inequality} and the statements thereafter.\\

Let $p\in P$ and $p'\in P$ be the points corresponding to the rays $\R\cdot  \iota_\cC(\mf{m})$ and $\R\cdot  \iota_\cC(\mf{m}')$, respectively.\\

For every $r,r'>0$, the multicurve $\mf{m}^r$ is a union of $r$ parallel copies of $\mf{m}$ (which can be drawn sufficiently close to $\mf{m}$), and the multicurve $(\mf{m}')^{r'}$ is a union of $r'$ parallel copies of $\mf{m}'$. Therefore, every intersection point of $\mf{m}$ and $\mf{m}'$ corresponds to a cluster of $rr'$ intersection points of $\mf{m}^r$ and $(\mf{m}')^{r'}$. Consider the smoothing $\mf{m}_{r,r',\#}$ of $(\mf{m}')^{r'}\cup (\mf{m})^{r}$, where the $rr'$ intersection points in each cluster are smoothed in the same way that $\mf{m}\cup \mf{m}'$ is smoothed to  $m_\#$ at the corresponding point. Any innermost bigon of $\mf{m}_{r,r',\#}$ with some curve $c\in \cC$ corresponds to an innermost bigon of  $m_\#$ with the same curve $c$. However, by \eqref{mforlemma}, the latter has no such bigon. Therefore,  
$$
\iota_\cC(\mf{m}_{r,r',\#})= \iota_\cC(\mf{m}^r)+\iota_\cC((\mf{m}')^{r'})=r \iota_\cC(\mf{m})+r'\iota_\cC(\mf{m}').
$$
In other words, the point $p_{r,r'}$ corresponding to the ray $\R\cdot \iota_\cC(\mf{m}_{r,r',\#})$ in $P\subset \De_{9g-10}$ satisfies
$$
p_{r,r'}=\frac{r}{r+r'} p + \frac{r'}{r+r'}p'.
$$
Since $r$ and $r'$ were arbitrary and the image of $\cS'$ is dense in $P$, we conclude that the line segment connecting $p$ and $p'$ is included in $P$; i.e., $p$ and $p'$ belong to the same cone in $\La$. That is a contradiction. \qed\medskip

This finishes the proof of Theorem~\ref{ThmC}.\qed\\

\section{Example of genus $2$}\label{sec:g2}
In this section, we dive into the details of the genus two example corresponding to the pants decomposition shown in Figure~\ref{example_fig}.\medskip

\begin{figure}
\begin{pspicture}(-7,-1)(0,1.5)
\psset{unit=1cm}

\psarc(0,0){1}{90}{270}
\psline(0,1)(2,1)
\psline(0,-1)(2,-1)
\psarc(2,0){1}{-90}{90}

\psarc(-0.5,0){.5}{-45}{45}   \psarc(2.5,0){.5}{135}{225}
\psarc(.3,0){.5}{145}{215} \psarc(1.7,0){.5}{-35}{35}

\psellipticarc[linecolor=red, linestyle=dashed, dash=1pt](-0.6,0)(0.4,0.1){0}{180} \rput(-0.6,.2){\small $a_1$}
\psellipticarc[linecolor=red](-0.6,0)(0.4,0.1){-180}{0}

\psellipticarc[linecolor=red, linestyle=dashed, dash=1pt](2.6,0)(0.4,0.1){0}{180} \rput(2.6,.2){\small $a_3$}
\psellipticarc[linecolor=red](2.6,0)(0.4,0.1){-180}{0}

\psellipticarc[linecolor=red, linestyle=dashed, dash=1pt](1,0)(1,0.1){0}{180} \rput(1,0.2){\small $a_2$}
\psellipticarc[linecolor=red](1,0)(1,0.1){180}{360}

\end{pspicture}
\caption{A pants decomposition of a genus two surface.}
\label{example_fig}
\end{figure} 

Let $F^\pm$ denote the upper and lower pair of pants. The DT coordinates are denoted by
$$
(n,t)=(n_1,n_2,n_3,t_1,t_2,t_3).
$$
The angle coordinates corresponding to each $a_i$ in both $F^+$ and $F^-$ coincide. Thus we denote the angle/twist coordinates by
$$
(\theta,t)=(\theta_1,\theta_2,\theta_3,t_1,t_2,t_3),
$$
where
$$
n_i=\theta_j+\theta_k \qquad \forall~i=1,2,3,\quad \{j,k\}=\{1,2,3\}\setminus\{i\}.
$$
We also have the $9g-9$ intersection numbers
$$
(n,n',n'')=(n_1,n_2,n_3,n'_1,n'_2,n'_3,n''_1,n''_2,n''_3),
$$
associated with the collection $\mathcal C$ in Theorem~\ref{ConeTheorem}.\medskip

The formulas for $n_i'$ and $n_i''$ in \eqref{FI} read
$$
\aligned
n_i'&=2\max \Big\{ \abs{t_i-\theta_i+C_i}, \theta_i-C_i\Big\}+4C_i-2\max\big\{\theta_i,0\big\},\\
n_i''&=2\max \Big\{ \abs{t_i-\theta_i-\theta_j-\theta_k+C_i}, \theta_i-C_i\Big\}+4C_i-2\max\big\{\theta_i,0\big\},
\endaligned
$$
where
$$
C_i=\max\big\{-\theta_j,0\big\}+\max\big\{-\theta_k,0\big\}+\max\big\{\theta_i,0\big\};
$$
see Remark~\ref{hf-in-thea}.
Note that
$$
\theta_i-C_i\leq 0\qquad \forall~i=1,2,3.
$$
Therefore, the formulas simplify to
$$
\aligned
n_i'&=2\abs{t_i-\theta_i+C_i}+4C_i-2\max\big\{\theta_i,0\big\},\\
n_i''&=2\abs{t_i-\theta_i-\theta_j-\theta_k+C_i}+4C_i-2\max\big\{\theta_i,0\big\}.
\endaligned
$$

There are four possibilities for the $\theta$-coordinates:
\begin{itemize}
\item Case $\Delta$: $n$ satisfies the triangle inequality and all $\theta_i$ are non-negative;
\item Cases $\Delta_i$: $n_i\geq n_j+n_k$ for some $i=1,2,3$, so that $\theta_i\leq 0$ and $\theta_j,\theta_k\geq 0$, subject to
$$
\theta_i+\theta_j\geq 0,\qquad \theta_i+\theta_k\geq 0,
$$
where $\{j,k\}=\{1,2,3\}\setminus\{i\}$.
\end{itemize}

In case $\Delta$, the formulas for $n'$ and $n''$ further simplify to
$$
\aligned
n_i'&=2 \abs{t_i}+2\theta_i,\\
n_i''&=2 \abs{t_i-(\theta_j+\theta_k)}+2\theta_i,
\endaligned
\qquad \forall~i=1,2,3.
$$
For each $i$, there are three domains of linearity depending on whether
$$
(A)\ t_i\leq 0,\qquad
(B)\ 0\leq t_i\leq \theta_j+\theta_k,\qquad
(C)\ t_i\geq \theta_j+\theta_k.
$$
Hence there are $27$ cones of type $\Delta$, denoted by $\Delta AAA,\ldots,\Delta CCC$.\medskip

For example, the generators of the (simplicial) cone $\Delta CCC$ are the vectors $v=(n,n',n'')$
$$
\aligned
&(\theta,t)=(0,0,0,1,0,0) \Rightarrow v_1=(0,0,0,2,0,0,2,0,0),\\
&(\theta,t)=(0,0,0,0,1,0)\Rightarrow v_2=(0,0,0,0,2,0,0,2,0),\\
&(\theta,t)=(0,0,0,0,0,1)\Rightarrow v_3=(0,0,0,0,0,2,0,0,2),\\
&(\theta,t)=(1,0,0,0,1,1)\Rightarrow v_4=(0,1,1,2,2,2,2,0,0),\\
&(\theta,t)=(0,1,0,1,0,1)\Rightarrow v_5=(1,0,1,2,2,2,0,2,0),\\
&(\theta,t)=(0,0,1,1,1,0)\Rightarrow v_6=(1,1,0,2,2,2,0,0,2).
\endaligned
$$
The degree of each term is the sum of the coordinates of $n$, $n'$, and $n''$:
$$
|v_1|=|v_2|=|v_3|=4,\qquad |v_4|=|v_5|=|v_6|=10.
$$
Thus the corresponding divisor component is the weighted projective space
$\P(4,4,4,10,10,10)$.\\

For $\De BBB$, on the other hand, the generators are the twelve degree $10$ vectors 
$$
\aligned
&(\theta,t)=(1,0,0,0,0,0)\Rightarrow v_1=(0,1,1,2,0,0,2,2,2)\\
&(\theta,t)=(1,0,0,0,1,0) \Rightarrow v_2=(0,1,1,2,2,0,2,0,2)\\
&(\theta,t)=(1,0,0,0,0,1) \Rightarrow v_3=(0,1,1,2,0,2,2,2,0)\\
&(\theta,t)=(1,0,0,0,1,1) \Rightarrow v_4=(0,1,1,2,2,2,2,0,0)\\
&(\theta,t)=(0,1,0,0,0,0)\Rightarrow u_1=(1,0,1,0,2,0,2,2,2)\\
&(\theta,t)=(0,1,0,1,0,0) \Rightarrow u_2=(1,0,1,2,2,0,0,2,2)\\
&(\theta,t)=(0,1,0,0,0,1) \Rightarrow u_3=(1,0,1,0,2,2,2,2,0)\\
&(\theta,t)=(0,1,0,1,0,1) \Rightarrow u_4=(1,0,1,2,2,2,0,2,0)\\
&(\theta,t)=(0,0,1,0,0,0)\Rightarrow w_1=(1,1,0,0,0,2,2,2,2)\\
&(\theta,t)=(0,0,1,1,0,0) \Rightarrow w_2=(1,1,0,2,0,2,0,2,2)\\
&(\theta,t)=(0,0,1,0,1,0) \Rightarrow w_3=(1,1,0,0,2,2,2,0,2)\\
&(\theta,t)=(0,0,1,1,1,0) \Rightarrow w_4=(1,1,0,2,2,2,0,0,2).
\endaligned
$$
A generating set of linear relations among these vectors is:
\[
\begin{aligned}
v_1 + v_4 &= v_2+v_3,\\
u_1 + u_4 &= u_2+u_3,\\
w_1 + w_4 &= w_2+w_3,\\
v_2 + u_1 &= v_1+u_2,\\
v_3 + w_1 &= v_1+w_3,\\
u_2 + w_1 &= u_1+w_2.
\end{aligned}
\]
Therefore, the corresponding boundary divisor is a toric complete intersection in the weighted projective space
$$\P^{11}\cong\P(10,\ldots,10).$$
In general, each cone of type $\De$ will have $6$ to $12$ generators, depending on the chosen intervals for $(t_1,t_2,t_3)$.\medskip

In case $\De_i$, the formulas for $n',n''$ read
$$
\aligned
n_i'&=2 \abs{t_i-\theta_i}\\
n_j'&=2\abs{t_j-\theta_i}+2\theta_j-4\theta_i\\
n_k'&=2\abs{t_k-\theta_i} +2\theta_k-4\theta_i\\
n_i''&=2 \abs{t_i-(\theta_i+\theta_j+\theta_k)}\\
n_j''&=2 \abs{t_j-2\theta_i-\theta_k}+2\theta_j-4\theta_i\\
n_k''&=2\abs{t_k-2\theta_i-\theta_j}+2\theta_k-4\theta_i\\
\endaligned
$$
and we get cones of various types as before. In total, there are $21$ extremal rays of degrees $4$, $10$, $38$, and $42$. Therefore, the whole boundary divisor $D_2=\partial \ov\cX_2$ sits in the weighted projective space 
$$
\P(4,4,4, 38, 38, 38, 42, 42, 42, 10, \cdots, 10).
$$    
\medskip

We conclude the paper by listing the equations describing the cones arising from the pants decomposition shown in Figure~\ref{otherexample_fig}.\medskip

\begin{figure}
\begin{pspicture}(-13.5,-1)(0,1.5)
\psset{unit=1cm}

\psarc(-7,0){1}{90}{270}
\psline(-7,1)(-5,1)
\psline(-7,-1)(-5,-1)
\psarc(-5,0){1}{-90}{90}

\psarc(-7.5,0){.5}{-45}{45}   
\psarc(-4.5,0){.5}{135}{225}
\psarc(-6.7,0){.5}{145}{215} 
\psarc(-5.3,0){.5}{-35}{35}

\psellipticarc[linecolor=red, linestyle=dashed, dash=1pt](-7.6,0)(0.4,0.1){0}{180} 
\rput(-7.6,.2){\small $a_L$}
\psellipticarc[linecolor=red](-7.6,0)(0.4,0.1){-180}{0}

\psellipticarc[linecolor=red, linestyle=dashed, dash=1pt](-4.4,0)(0.4,0.1){0}{180} 
\rput(-4.4,.2){\small $a_R$}
\psellipticarc[linecolor=red](-4.4,0)(0.4,0.1){-180}{0}

\psellipticarc[linecolor=red](-6,0)(0.2,1){90}{270} 
\rput(-6.4,0){\small $a$}
\psellipticarc[linecolor=red, linestyle=dashed, dash=1pt](-6,0)(0.2,1){-90}{90}

\end{pspicture}
\caption{A pants decomposition of a genus two surface with two folded pants.}
\label{otherexample_fig}
\end{figure}

Let $F_L$ and $F_R$ denote the left and right pairs of pants, respectively. The DT coordinates are denoted by
$$
(n,t)=(n,n_L,n_R,t,t_L,t_R).
$$
The (only) angle coordinates corresponding to $a_L$ and $a_R$ coincide and are given by
$$
\theta=\frac{n}{2}.
$$
The angle coordinates corresponding to $(F_L,a)$ and $(F_R,a)$ are
$$
\theta_L=n_L-\theta \quad\text{and}\quad \theta_R=n_R-\theta.
$$
Thus the tuple of angle/twist coordinates is
$$
(\theta,\theta_L,\theta_R,t,t_L,t_R),
$$
subject to the relations
$$
\aligned
&\theta\geq 0,\qquad \theta+\theta_L\geq 0,\qquad \theta+\theta_R\geq 0,\\
&\theta=0 \Rightarrow t\geq 0 \ \text{or}\ t\leq 0 \quad (\text{depending on the chosen cone}),\\
&\theta+\theta_L=0 \Rightarrow t_L\geq 0 \ \text{or}\ t_L\leq 0 \quad (\text{depending on the chosen cone}),\\
&\theta+\theta_R=0 \Rightarrow t_R\geq 0 \ \text{or}\ t_R\leq 0 \quad (\text{depending on the chosen cone}).
\endaligned
$$
We also have the $9g-9$ intersection numbers
$$
(n,n',n'')=(n,n_L,n_R,n',n'_L,n'_R,n'',n''_L,n''_R),
$$
associated with the collection $\mathcal C$ in Theorem~\ref{ConeTheorem}.\\

The formulas for $n'$ and $n''$ in \eqref{FI} and \eqref{FIfolded} yield
$$
\aligned
n'_L&=\abs{t_L}+\max\big\{-\theta_L,0\big\},\\
n'_R&=\abs{t_R}+\max\big\{-\theta_R,0\big\},\\
n''_L&=\abs{t_L-(\theta+\theta_L)}+\max\big\{-\theta_L,0\big\},\\
n''_R&=\abs{t_R-(\theta+\theta_R)}+\max\big\{-\theta_R,0\big\},\\
n'&=\abs{2t-2C}+\max\big\{\theta_L,0\big\}+\max\big\{\theta_R,0\big\},\\
n''&=\abs{2t-4\theta-2C}+\max\big\{\theta_L,0\big\}+\max\big\{\theta_R,0\big\},
\endaligned
$$
where
$$
C=\theta_L+\theta_R-\max\big\{\theta_L,0\big\}-\max\big\{\theta_R,0\big\}\leq 0.
$$

The max functions depend on four possible scenarios, according to whether the triangle inequalities
$2n_L\geq n$ and $2n_R\geq n$ hold. In each case, there are three regions for each of the variables
$t$, $t_L$, and $t_R$. Altogether, this yields $4\times 27$ cones of various combinatorial types as before.

\bibliographystyle{plain}

\bibliographystyle{amsalpha}

\end{document}